\catcode`@=11
\def\@height{height}
\def\@depth{depth}
\def\@width{width}

\newcount\@tempcnta
\newcount\@tempcntb

\newdimen\@tempdima
\newdimen\@tempdimb

\newbox\@tempboxa

\def\@ifnextchar#1#2#3{\let\@tempe #1\def\@tempa{#2}\def\@tempb{#3}\futurelet
    \@tempc\@ifnch}
\def\@ifnch{\ifx \@tempc \@sptoken \let\@tempd\@xifnch
      \else \ifx \@tempc \@tempe\let\@tempd\@tempa\else\let\@tempd\@tempb\fi
      \fi \@tempd}
\def\@ifstar#1#2{\@ifnextchar *{\def\@tempa*{#1}\@tempa}{#2}}

\def\@whilenoop#1{}
\def\@whilenum#1\do #2{\ifnum #1\relax #2\relax\@iwhilenum{#1\relax 
     #2\relax}\fi}
\def\@iwhilenum#1{\ifnum #1\let\@nextwhile=\@iwhilenum 
         \else\let\@nextwhile=\@whilenoop\fi\@nextwhile{#1}}

\def\@whiledim#1\do #2{\ifdim #1\relax#2\@iwhiledim{#1\relax#2}\fi}
\def\@iwhiledim#1{\ifdim #1\let\@nextwhile=\@iwhiledim 
        \else\let\@nextwhile=\@whilenoop\fi\@nextwhile{#1}}

\newdimen\@wholewidth
\newdimen\@halfwidth
\newdimen\unitlength \unitlength =1pt
\newbox\@picbox
\newdimen\@picht

\def\@nnil{\@nil}
\def\@empty{}
\def\@fornoop#1\@@#2#3{}

\def\@for#1:=#2\do#3{\edef\@fortmp{#2}\ifx\@fortmp\@empty \else
    \expandafter\@forloop#2,\@nil,\@nil\@@#1{#3}\fi}

\def\@forloop#1,#2,#3\@@#4#5{\def#4{#1}\ifx #4\@nnil \else
       #5\def#4{#2}\ifx #4\@nnil \else#5\@iforloop #3\@@#4{#5}\fi\fi}

\def\@iforloop#1,#2\@@#3#4{\def#3{#1}\ifx #3\@nnil 
       \let\@nextwhile=\@fornoop \else
      #4\relax\let\@nextwhile=\@iforloop\fi\@nextwhile#2\@@#3{#4}}

\def\@tfor#1:=#2\do#3{\xdef\@fortmp{#2}\ifx\@fortmp\@empty \else
    \@tforloop#2\@nil\@nil\@@#1{#3}\fi}
\def\@tforloop#1#2\@@#3#4{\def#3{#1}\ifx #3\@nnil 
       \let\@nextwhile=\@fornoop \else
      #4\relax\let\@nextwhile=\@tforloop\fi\@nextwhile#2\@@#3{#4}}

\def\@makepicbox(#1,#2){\leavevmode\@ifnextchar 
   [{\@imakepicbox(#1,#2)}{\@imakepicbox(#1,#2)[]}}

\long\def\@imakepicbox(#1,#2)[#3]#4{\vbox to#2\unitlength
   {\let\mb@b\vss \let\mb@l\hss\let\mb@r\hss
    \let\mb@t\vss
    \@tfor\@tempa :=#3\do{\expandafter\let
        \csname mb@\@tempa\endcsname\relax}%
\mb@t\hbox to #1\unitlength{\mb@l #4\mb@r}\mb@b}}

\def\picture(#1,#2){\@ifnextchar({\@picture(#1,#2)}{\@picture(#1,#2)(0,0)}}

\def\@picture(#1,#2)(#3,#4){\@picht #2\unitlength
\setbox\@picbox\hbox to #1\unitlength\bgroup 
\hskip -#3\unitlength \lower #4\unitlength \hbox\bgroup\ignorespaces}

\def\endpicture{\egroup\hss\egroup\ht\@picbox\@picht
\dp\@picbox\z@\leavevmode\box\@picbox}

\long\def\put(#1,#2)#3{\@killglue\raise#2\unitlength\hbox to \z@{\kern
#1\unitlength #3\hss}\ignorespaces}

\long\def\multiput(#1,#2)(#3,#4)#5#6{\@killglue\@multicnt=#5\relax
\@xdim=#1\unitlength
\@ydim=#2\unitlength
\@whilenum \@multicnt > 0\do
{\raise\@ydim\hbox to \z@{\kern
\@xdim #6\hss}\advance\@multicnt \m@ne\advance\@xdim
#3\unitlength\advance\@ydim #4\unitlength}\ignorespaces}

\def\@killglue{\unskip\@whiledim \lastskip >\z@\do{\unskip}}

\def\thinlines{\let\@linefnt\tenln \let\@circlefnt\tencirc
  \@wholewidth\fontdimen8\tenln \@halfwidth .5\@wholewidth}
\def\thicklines{\let\@linefnt\tenlnw \let\@circlefnt\tencircw
  \@wholewidth\fontdimen8\tenlnw \@halfwidth .5\@wholewidth}

\def\linethickness#1{\@wholewidth #1\relax \@halfwidth .5\@wholewidth}

\def\shortstack{\@ifnextchar[{\@shortstack}{\@shortstack[c]}}

\def\@shortstack[#1]{\leavevmode
\vbox\bgroup\baselineskip-1pt\lineskip 3pt\let\mb@l\hss
\let\mb@r\hss \expandafter\let\csname mb@#1\endcsname\relax
\let\\\@stackcr\@ishortstack}

\def\@ishortstack#1{\halign{\mb@l ##\unskip\mb@r\cr #1\crcr}\egroup}

\def\@stackcr{\@ifstar{\@ixstackcr}{\@ixstackcr}}
\def\@ixstackcr{\@ifnextchar[{\@istackcr}{\cr\ignorespaces}}

\def\@istackcr[#1]{\cr\noalign{\vskip #1}\ignorespaces}

\newif\if@negarg

\def\droite(#1,#2)#3{\@xarg #1\relax \@yarg #2\relax
\@linelen=#3\unitlength
\ifnum\@xarg =0 \@vline 
  \else \ifnum\@yarg =0 \@hline \else \@sline\fi
\fi}

\def\@sline{\ifnum\@xarg< 0 \@negargtrue \@xarg -\@xarg \@yyarg -\@yarg
  \else \@negargfalse \@yyarg \@yarg \fi
\ifnum \@yyarg >0 \@tempcnta\@yyarg \else \@tempcnta -\@yyarg \fi
\ifnum\@tempcnta>6 \@badlinearg\@tempcnta0 \fi
\ifnum\@xarg>6 \@badlinearg\@xarg 1 \fi
\setbox\@linechar\hbox{\@linefnt\@getlinechar(\@xarg,\@yyarg)}%
\ifnum \@yarg >0 \let\@upordown\raise \@clnht\z@
   \else\let\@upordown\lower \@clnht \ht\@linechar\fi
\@clnwd=\wd\@linechar
\if@negarg \hskip -\wd\@linechar \def\@tempa{\hskip -2\wd\@linechar}\else
     \let\@tempa\relax \fi
\@whiledim \@clnwd <\@linelen \do
  {\@upordown\@clnht\copy\@linechar
   \@tempa
   \advance\@clnht \ht\@linechar
   \advance\@clnwd \wd\@linechar}%
\advance\@clnht -\ht\@linechar
\advance\@clnwd -\wd\@linechar
\@tempdima\@linelen\advance\@tempdima -\@clnwd
\@tempdimb\@tempdima\advance\@tempdimb -\wd\@linechar
\if@negarg \hskip -\@tempdimb \else \hskip \@tempdimb \fi
\multiply\@tempdima \@m
\@tempcnta \@tempdima \@tempdima \wd\@linechar \divide\@tempcnta \@tempdima
\@tempdima \ht\@linechar \multiply\@tempdima \@tempcnta
\divide\@tempdima \@m
\advance\@clnht \@tempdima
\ifdim \@linelen <\wd\@linechar
   \hskip \wd\@linechar
  \else\@upordown\@clnht\copy\@linechar\fi}

\def\@hline{\ifnum \@xarg <0 \hskip -\@linelen \fi
\vrule \@height \@halfwidth \@depth \@halfwidth \@width \@linelen
\ifnum \@xarg <0 \hskip -\@linelen \fi}

\def\@getlinechar(#1,#2){\@tempcnta#1\relax\multiply\@tempcnta 8
\advance\@tempcnta -9 \ifnum #2>0 \advance\@tempcnta #2\relax\else
\advance\@tempcnta -#2\relax\advance\@tempcnta 64 \fi
\char\@tempcnta}

\def\vector(#1,#2)#3{\@xarg #1\relax \@yarg #2\relax
\@tempcnta \ifnum\@xarg<0 -\@xarg\else\@xarg\fi
\ifnum\@tempcnta<5\relax
\@linelen=#3\unitlength
\ifnum\@xarg =0 \@vvector 
  \else \ifnum\@yarg =0 \@hvector \else \@svector\fi
\fi
\else\@badlinearg\fi}

\def\@hvector{\@hline\hbox to 0pt{\@linefnt 
\ifnum \@xarg <0 \@getlarrow(1,0)\hss\else
    \hss\@getrarrow(1,0)\fi}}

\def\@vvector{\ifnum \@yarg <0 \@downvector \else \@upvector \fi}

\def\@svector{\@sline
\@tempcnta\@yarg \ifnum\@tempcnta <0 \@tempcnta=-\@tempcnta\fi
\ifnum\@tempcnta <5
  \hskip -\wd\@linechar
  \@upordown\@clnht \hbox{\@linefnt  \if@negarg 
  \@getlarrow(\@xarg,\@yyarg) \else \@getrarrow(\@xarg,\@yyarg) \fi}%
\else\@badlinearg\fi}

\def\@getlarrow(#1,#2){\ifnum #2 =\z@ \@tempcnta='33\else
\@tempcnta=#1\relax\multiply\@tempcnta \sixt@@n \advance\@tempcnta
-9 \@tempcntb=#2\relax\multiply\@tempcntb \tw@
\ifnum \@tempcntb >0 \advance\@tempcnta \@tempcntb\relax
\else\advance\@tempcnta -\@tempcntb\advance\@tempcnta 64
\fi\fi\char\@tempcnta}

\def\@getrarrow(#1,#2){\@tempcntb=#2\relax
\ifnum\@tempcntb < 0 \@tempcntb=-\@tempcntb\relax\fi
\ifcase \@tempcntb\relax \@tempcnta='55 \or 
\ifnum #1<3 \@tempcnta=#1\relax\multiply\@tempcnta
24 \advance\@tempcnta -6 \else \ifnum #1=3 \@tempcnta=49
\else\@tempcnta=58 \fi\fi\or 
\ifnum #1<3 \@tempcnta=#1\relax\multiply\@tempcnta
24 \advance\@tempcnta -3 \else \@tempcnta=51\fi\or 
\@tempcnta=#1\relax\multiply\@tempcnta
\sixt@@n \advance\@tempcnta -\tw@ \else
\@tempcnta=#1\relax\multiply\@tempcnta
\sixt@@n \advance\@tempcnta 7 \fi\ifnum #2<0 \advance\@tempcnta 64 \fi
\char\@tempcnta}

\def\@vline{\ifnum \@yarg <0 \@downline \else \@upline\fi}

\def\@upline{\hbox to \z@{\hskip -\@halfwidth \vrule \@width \@wholewidth
   \@height \@linelen \@depth \z@\hss}}

\def\@downline{\hbox to \z@{\hskip -\@halfwidth \vrule \@width \@wholewidth
   \@height \z@ \@depth \@linelen \hss}}

\def\@upvector{\@upline\setbox\@tempboxa\hbox{\@linefnt\char'66}\raise 
     \@linelen \hbox to\z@{\lower \ht\@tempboxa\box\@tempboxa\hss}}

\def\@downvector{\@downline\lower \@linelen
      \hbox to \z@{\@linefnt\char'77\hss}}

\def\dashbox#1(#2,#3){\leavevmode\hbox to \z@{\baselineskip \z@%
\lineskip \z@%
\@dashdim=#2\unitlength%
\@dashcnt=\@dashdim \advance\@dashcnt 200
\@dashdim=#1\unitlength\divide\@dashcnt \@dashdim
\ifodd\@dashcnt\@dashdim=\z@%
\advance\@dashcnt \@ne \divide\@dashcnt \tw@ 
\else \divide\@dashdim \tw@ \divide\@dashcnt \tw@
\advance\@dashcnt \m@ne
\setbox\@dashbox=\hbox{\vrule \@height \@halfwidth \@depth \@halfwidth
\@width \@dashdim}\put(0,0){\copy\@dashbox}%
\put(0,#3){\copy\@dashbox}%
\put(#2,0){\hskip-\@dashdim\copy\@dashbox}%
\put(#2,#3){\hskip-\@dashdim\box\@dashbox}%
\multiply\@dashdim 3 
\fi
\setbox\@dashbox=\hbox{\vrule \@height \@halfwidth \@depth \@halfwidth
\@width #1\unitlength\hskip #1\unitlength}\@tempcnta=0
\put(0,0){\hskip\@dashdim \@whilenum \@tempcnta <\@dashcnt
\do{\copy\@dashbox\advance\@tempcnta \@ne }}\@tempcnta=0
\put(0,#3){\hskip\@dashdim \@whilenum \@tempcnta <\@dashcnt
\do{\copy\@dashbox\advance\@tempcnta \@ne }}%
\@dashdim=#3\unitlength%
\@dashcnt=\@dashdim \advance\@dashcnt 200
\@dashdim=#1\unitlength\divide\@dashcnt \@dashdim
\ifodd\@dashcnt \@dashdim=\z@%
\advance\@dashcnt \@ne \divide\@dashcnt \tw@
\else
\divide\@dashdim \tw@ \divide\@dashcnt \tw@
\advance\@dashcnt \m@ne
\setbox\@dashbox\hbox{\hskip -\@halfwidth
\vrule \@width \@wholewidth 
\@height \@dashdim}\put(0,0){\copy\@dashbox}%
\put(#2,0){\copy\@dashbox}%
\put(0,#3){\lower\@dashdim\copy\@dashbox}%
\put(#2,#3){\lower\@dashdim\copy\@dashbox}%
\multiply\@dashdim 3
\fi
\setbox\@dashbox\hbox{\vrule \@width \@wholewidth 
\@height #1\unitlength}\@tempcnta0
\put(0,0){\hskip -\@halfwidth \vbox{\@whilenum \@tempcnta < \@dashcnt
\do{\vskip #1\unitlength\copy\@dashbox\advance\@tempcnta \@ne }%
\vskip\@dashdim}}\@tempcnta0
\put(#2,0){\hskip -\@halfwidth \vbox{\@whilenum \@tempcnta< \@dashcnt
\relax\do{\vskip #1\unitlength\copy\@dashbox\advance\@tempcnta \@ne }%
\vskip\@dashdim}}}\@makepicbox(#2,#3)}

\newif\if@ovt 
\newif\if@ovb 
\newif\if@ovl 
\newif\if@ovr 
\newdimen\@ovxx
\newdimen\@ovyy
\newdimen\@ovdx
\newdimen\@ovdy
\newdimen\@ovro
\newdimen\@ovri

\def\@getcirc#1{\@tempdima #1\relax \advance\@tempdima 2pt\relax
  \@tempcnta\@tempdima
  \@tempdima 4pt\relax \divide\@tempcnta\@tempdima
  \ifnum \@tempcnta > 10\relax \@tempcnta 10\relax\fi
  \ifnum \@tempcnta >\z@ \advance\@tempcnta\m@ne
    \else \@warning{Oval too small}\fi
  \multiply\@tempcnta 4\relax
  \setbox \@tempboxa \hbox{\@circlefnt
  \char \@tempcnta}\@tempdima \wd \@tempboxa}

\def\@put#1#2#3{\raise #2\hbox to \z@{\hskip #1#3\hss}}

\def\oval(#1,#2){\@ifnextchar[{\@oval(#1,#2)}{\@oval(#1,#2)[]}}

\def\@oval(#1,#2)[#3]{\begingroup\boxmaxdepth \maxdimen
  \@ovttrue \@ovbtrue \@ovltrue \@ovrtrue
  \@tfor\@tempa :=#3\do{\csname @ov\@tempa false\endcsname}\@ovxx
  #1\unitlength \@ovyy #2\unitlength
  \@tempdimb \ifdim \@ovyy >\@ovxx \@ovxx\else \@ovyy \fi
  \advance \@tempdimb -2pt\relax  %%%% added 7 Dec 89
  \@getcirc \@tempdimb
  \@ovro \ht\@tempboxa \@ovri \dp\@tempboxa
  \@ovdx\@ovxx \advance\@ovdx -\@tempdima \divide\@ovdx \tw@
  \@ovdy\@ovyy \advance\@ovdy -\@tempdima \divide\@ovdy \tw@
  \@circlefnt \setbox\@tempboxa
  \hbox{\if@ovr \@ovvert32\kern -\@tempdima \fi
  \if@ovl \kern \@ovxx \@ovvert01\kern -\@tempdima \kern -\@ovxx \fi
  \if@ovt \@ovhorz \kern -\@ovxx \fi
  \if@ovb \raise \@ovyy \@ovhorz \fi}\advance\@ovdx\@ovro
  \advance\@ovdy\@ovro \ht\@tempboxa\z@ \dp\@tempboxa\z@
  \@put{-\@ovdx}{-\@ovdy}{\box\@tempboxa}%
  \endgroup}

\def\@ovvert#1#2{\vbox to \@ovyy{%
    \if@ovb \@tempcntb \@tempcnta \advance \@tempcntb by #1\relax
      \kern -\@ovro \hbox{\char \@tempcntb}\nointerlineskip
    \else \kern \@ovri \kern \@ovdy \fi
    \leaders\vrule width \@wholewidth\vfil \nointerlineskip
    \if@ovt \@tempcntb \@tempcnta \advance \@tempcntb by #2\relax
      \hbox{\char \@tempcntb}%
    \else \kern \@ovdy \kern \@ovro \fi}}

\def\@ovhorz{\hbox to \@ovxx{\kern \@ovro
    \if@ovr \else \kern \@ovdx \fi
    \leaders \hrule height \@wholewidth \hfil
    \if@ovl \else \kern \@ovdx \fi
    \kern \@ovri}}

\def\circle{\@ifstar{\@dot}{\@circle}}
\def\@circle#1{\begingroup \boxmaxdepth \maxdimen \@tempdimb #1\unitlength
   \ifdim \@tempdimb >15.5pt\relax \@getcirc\@tempdimb
      \@ovro\ht\@tempboxa 
     \setbox\@tempboxa\hbox{\@circlefnt
      \advance\@tempcnta\tw@ \char \@tempcnta
      \advance\@tempcnta\m@ne \char \@tempcnta \kern -2\@tempdima
      \advance\@tempcnta\tw@
      \raise \@tempdima \hbox{\char\@tempcnta}\raise \@tempdima
        \box\@tempboxa}\ht\@tempboxa\z@ \dp\@tempboxa\z@
      \@put{-\@ovro}{-\@ovro}{\box\@tempboxa}%
   \else  \@circ\@tempdimb{96}\fi\endgroup}

\def\@dot#1{\@tempdimb #1\unitlength \@circ\@tempdimb{112}}

\def\@circ#1#2{\@tempdima #1\relax \advance\@tempdima .5pt\relax
   \@tempcnta\@tempdima \@tempdima 1pt\relax
   \divide\@tempcnta\@tempdima 
   \ifnum\@tempcnta > 15\relax \@tempcnta 15\relax \fi    
   \ifnum \@tempcnta >\z@ \advance\@tempcnta\m@ne\fi
   \advance\@tempcnta #2\relax
   \@circlefnt \char\@tempcnta}

%INITIALIZATION
\font\tenln line10
%\font\tencirc circle10
\font\tencirc lcircle10
\font\tenlnw linew10
%\font\tencircw circlew10
\font\tencircw lcirclew10

\thinlines   

\newcount\@xarg
\newcount\@yarg
\newcount\@yyarg
\newcount\@multicnt 
\newdimen\@xdim
\newdimen\@ydim
\newbox\@linechar
\newdimen\@linelen
\newdimen\@clnwd
\newdimen\@clnht
\newdimen\@dashdim
\newbox\@dashbox
\newcount\@dashcnt
\catcode`@=12
%%% Local Variables: 
%%% mode: plain-tex
%%% TeX-master: t
%%% End: 

%This is trialgebra.macros

\overfullrule=0pt
\magnification=1200
\hsize=11.25cm    
\vsize=19cm
\hoffset=1cm

\font\grand=cmr10 at 12pt

\def\N{\noindent}
\def\S{\smallskip \par}
\def\M{\medskip \par}
\def\B{\bigskip \par}
\def\BB{\bigskip \bigskip \par}
\def\P{\noindent {\it Proof.} }

\catcode`@=11
\font\@linefnt linew10 at 2.4pt
\catcode`@=12

\font\grand=cmr10 at 14pt

\def\u{\underline }
\def\P{{\cal P}}
\def\K{{\cal K}}

%opŽration gauche
\def\g{\mathop{\dashv}\nolimits}

%opŽration droite
\def\d{\mathop{\vdash}\nolimits}

%opŽration milieu
\def\m{\mathop{\perp}\nolimits}

%opŽration gauche dendriform
\def\l{\mathop{\prec}\nolimits}

%opŽration droite dendriform
\def\r{\mathop{\succ}\nolimits}

%puissance moins un

\def\N{\noindent}
\def\S{\smallskip \par}
\def\M{\medskip \par}
\def\B{\bigskip \par}
\def\BB{\bigskip \bigskip \par}
%\def\P{\noindent {\it Proof.} }

%lettres grecques
\def\aa{\alpha}
\def\bb{\beta}
\def\cc{\gamma}

\def\dd{\delta}
\def\DD{\Delta}

\def\ss{\sigma}
\def\SS{\Sigma}
\def\oo{\omega}

%carrŽ de fin de dŽmonstration
\def\sqr#1#2{{\vcenter{\vbox{\hrule height.#2pt
\hbox{\vrule width .#2pt height#1pt \kern#1pt
\vrule width.#2pt}
\hrule height.#2pt}}}}
\def\square{\mathchoice\sqr64\sqr64\sqr{4.2}3\sqr33}

\def \Id{\mathop{\rm Id}\nolimits}

\def \Vect{\mathop{\rm Vect}\nolimits}

%tenseur
\def\t{\otimes }

%V tenseur n
\def\tT#1#2{#1^{\otimes #2}}

%vecteur horizontal #1= variable #2= premier indice
%#3=dernier indice

%#3=dernier indice

\def\arbreA{\kern-0.4ex
\hbox{\unitlength=.25pt
\picture(60,40)(0,0)
\put(30,0){\droite(0,1){10}}
\put(30,10){\droite(-1,1){30}}
\put(30,10){\droite(1,1){30}}
\endpicture}\kern 0.4ex}

\def\arbreB{\kern-0.4ex
\hbox{\unitlength=.25pt
\picture(60,40)(0,0)
\put(30,0){\droite(0,1){10}}
\put(30,10){\droite(-1,1){30}}
\put(30,10){\droite(1,1){30}}
\put(15,25){\droite(1,1){15}}
\endpicture}\kern 0.4ex}

\def\arbreC{\kern-0.4ex
\hbox{\unitlength=.25pt
\picture(60,40)(0,0)
\put(30,0){\droite(0,1){10}}
\put(30,10){\droite(-1,1){30}}
\put(30,10){\droite(1,1){30}}
\put(45,25){\droite(-1,1){15}}
\endpicture}\kern 0.4ex}

\def\arbreBC{\kern-0.4ex
\hbox{\unitlength=.25pt
\picture(60,40)(0,0)
\put(30,0){\droite(0,1){40}}
\put(30,10){\droite(-1,1){30}}
\put(30,10){\droite(1,1){30}}
\endpicture}\kern 0.4ex}

\def\arbreun{\kern-0.4ex
\hbox{\unitlength=.25pt
\picture(60,40)(0,0)
\put(30,0){\droite(0,1){10}}
\put(30,10){\droite(-1,1){30}}
\put(30,10){\droite(1,1){30}}
\put(20,20){\droite(1,1){20}}
\put(10,30){\droite(1,1){10}}
\endpicture}\kern 0.4ex}

\def\arbredeux{\kern-0.4ex
\hbox{\unitlength=.25pt
\picture(60,40)(0,0)
\put(30,0){\droite(0,1){10}}
\put(30,10){\droite(-1,1){30}}
\put(30,10){\droite(1,1){30}}
\put(20,20){\droite(1,1){20}}
\put(30,30){\droite(-1,1){10}}
\endpicture}\kern 0.4ex}

\def\arbretrois{\kern-0.4ex
\hbox{\unitlength=.25pt
\picture(60,40)(0,0)
\put(30,0){\droite(0,1){10}}
\put(30,10){\droite(-1,1){30}}
\put(30,10){\droite(1,1){30}}
\put(50,30){\droite(-1,1){10}}
\put(10,30){\droite(1,1){10}}
\endpicture}\kern 0.4ex}

\def\arbrequatre{\kern-0.4ex
\hbox{\unitlength=.25pt
\picture(60,40)(0,0)
\put(30,0){\droite(0,1){10}}
\put(30,10){\droite(-1,1){30}}
\put(30,10){\droite(1,1){30}}
\put(40,20){\droite(-1,1){20}}
\put(30,30){\droite(1,1){10}}
\endpicture}\kern 0.4ex}

\def\arbrecinq{\kern-0.4ex
\hbox{\unitlength=.25pt
\picture(60,40)(0,0)
\put(30,0){\droite(0,1){10}}
\put(30,10){\droite(-1,1){30}}
\put(30,10){\droite(1,1){30}}
\put(40,20){\droite(-1,1){20}}
\put(50,30){\droite(-1,1){10}}
\endpicture}\kern 0.4ex}

\def\arbreuut{\kern-0.4ex
\hbox{\unitlength=.25pt
\picture(60,40)(0,0)
\put(30,0){\droite(0,1){10}}
\put(30,10){\droite(-1,1){30}}
\put(30,10){\droite(1,1){30}}
\put(20,20){\droite(0,1){20}}
\put(20,20){\droite(1,1){20}}
\endpicture}\kern 0.4ex}

\def\arbretut{\kern-0.4ex
\hbox{\unitlength=.25pt
\picture(60,40)(0,0)
\put(30,0){\droite(0,1){10}}
\put(30,10){\droite(-1,1){30}}
\put(30,10){\droite(1,1){30}}
\put(30,10){\droite(0,1){20}}
\put(30,30){\droite(1,1){10}}
\put(30,30){\droite(-1,1){10}}
\endpicture}\kern 0.4ex}

\def\arbretuu{\kern-0.4ex
\hbox{\unitlength=.25pt
\picture(60,40)(0,0)
\put(30,0){\droite(0,1){10}}
\put(30,10){\droite(-1,1){30}}
\put(30,10){\droite(1,1){30}}
\put(40,20){\droite(0,1){20}}
\put(40,20){\droite(-1,1){20}}
\endpicture}\kern 0.4ex}

\def\arbreutt{\kern-0.4ex
\hbox{\unitlength=.25pt
\picture(60,40)(0,0)
\put(30,0){\droite(0,1){10}}
\put(30,10){\droite(-1,1){30}}
\put(30,10){\droite(1,1){30}}
\put(10,30){\droite(1,1){10}}
\put(30,10){\droite(0,1){30}}
\endpicture}\kern 0.4ex}

\def\arbrettu{\kern-0.4ex
\hbox{\unitlength=.25pt
\picture(60,40)(0,0)
\put(30,0){\droite(0,1){10}}
\put(30,10){\droite(-1,1){30}}
\put(30,10){\droite(1,1){30}}
\put(50,30){\droite(-1,1){10}}
\put(30,10){\droite(0,1){30}}
\endpicture}\kern 0.4ex}

\def\arbrettt{\kern-0.4ex
\hbox{\unitlength=.25pt
\picture(60,40)(0,0)
\put(30,0){\droite(0,1){10}}
\put(30,10){\droite(-1,1){30}}
\put(30,10){\droite(1,1){30}}
\put(30,10){\droite(-1,2){15}}
\put(30,10){\droite(1,2){15}}
\endpicture}\kern 0.4ex}

\def\espaceXxxx{\kern-0.4ex
\hbox{\unitlength=.25pt
\picture(600,400)(0,0)
\put(0,0){$\bullet$}
\put(300,0){$\bullet$}
\put(600,0){$\bullet$}
\put(300,300){$\bullet$}
\put(0,10){\droite(1,0){600}}
\put(310,10){\droite(0,1){300}}
\endpicture}\kern 0.4ex}

\def\simplexes{\kern-0.4ex
\hbox{\unitlength=.25pt
\picture(800,250)(0,0)
\put(0,0){$\Delta^0$}
\put(10,50){$\bullet$}
\put(10,100){$0$}

\put(200,0){$\Delta^1$}
\put(160,50){$\bullet$}
\put(170,60){\droite(1,0){100}}
\put(260,50){$\bullet$}
\put(160,100){$0$}
\put(260,100){$1$}

\put(450,0){$\Delta^2$}
\put(410,50){$\bullet$}
\put(420,60){\droite(1,0){100}}
\put(420,60){\droite(1,2){50}}
\put(520,60){\droite(-1,2){50}}
\put(410,50){$\bullet$}
\put(390,80){$0$}
\put(460,150){$\bullet$}
\put(530,80){$1$}
\put(510,50){$\bullet$}
\put(460,180){$2$}

\endpicture}\kern 0.4ex}

\def\stasheff{\kern-0.4ex
\hbox{\unitlength=.25pt
\picture(800,270)(0,0)
\put(0,0){$\K^0$}
\put(10,50){$\bullet$}
\put(0,100){$\arbreA$}

\put(200,0){$\K^1$}
\put(160,50){$\bullet$}
\put(170,60){\droite(1,0){120}}
\put(280,50){$\bullet$}
\put(140,100){$\arbreB$}
\put(270,100){$\arbreC$}
\put(210,70){$\arbreBC$}

\put(590,0){$\K^2$}
\put(520,130){\droite(1,1){50}}
\put(570,80){\droite(3,1){50}}
\put(620,100){\droite(0,1){60}}
\put(570,180){\droite(3,-1){50}}
\put(520,130){\droite(1,-1){50}}

\put(610,150){$\bullet$}
\put(640,150){$\arbredeux$}

\put(560,170){$\bullet$}
\put(540,200){$\arbreun$}

\put(510,120){$\bullet$}
\put(460,110){$\arbretrois$}

\put(560,70){$\bullet$}
\put(640,70){$\arbrequatre$}

\put(610,90){$\bullet$}
\put(530,20){$\arbrecinq$}

\endpicture}\kern 0.4ex}

\def\figure5A{\kern-0.4ex
\hbox{\unitlength=.25pt
\picture(120,80)(0,0)
\put(60,0){\droite(0,1){10}}
\put(60,10){\droite(-1,1){60}}
\put(60,10){\droite(1,1){60}}
\put(70,20){\droite(-1,1){50}}
\put(45,50){$\cdots$}
\endpicture}\kern 0.4ex}

\def\arbrequdt{\kern-0.4ex
\hbox{\unitlength=.25pt
\picture(80,50)(0,0)
\put(40,0){\droite(0,1){10}}
\put(40,10){\droite(-1,1){40}}
\put(40,10){\droite(1,1){40}}
\put(50,20){\droite(-1,1){30}}
\put(40,30){\droite(1,1){20}}
\put(30,40){\droite(1,1){10}}
\endpicture}\kern 0.4ex}

\def\arbreqdut{\kern-0.4ex
\hbox{\unitlength=.25pt
\picture(80,50)(0,0)
\put(40,0){\droite(0,1){10}}
\put(40,10){\droite(-1,1){40}}
\put(40,10){\droite(1,1){40}}
\put(50,20){\droite(-1,1){30}}
\put(40,30){\droite(1,1){20}}
\put(50,40){\droite(-1,1){10}}
\endpicture}\kern 0.4ex}

\def\arbreqtud{\kern-0.4ex
\hbox{\unitlength=.25pt
\picture(80,50)(0,0)
\put(40,0){\droite(0,1){10}}
\put(40,10){\droite(-1,1){40}}
\put(40,10){\droite(1,1){40}}
\put(50,20){\droite(-1,1){30}}
\put(60,30){\droite(-1,1){20}}
\put(50,40){\droite(1,1){10}}
\endpicture}\kern 0.4ex}

\def\arbreqtdu{\kern-0.4ex
\hbox{\unitlength=.25pt
\picture(80,50)(0,0)
\put(40,0){\droite(0,1){10}}
\put(40,10){\droite(-1,1){40}}
\put(40,10){\droite(1,1){40}}
\put(50,20){\droite(-1,1){30}}
\put(60,30){\droite(-1,1){20}}
\put(70,40){\droite(-1,1){10}}
\endpicture}\kern 0.4ex}

\def\arbrequtu{\kern-0.4ex
\hbox{\unitlength=.25pt
\picture(80,50)(0,0)
\put(40,0){\droite(0,1){10}}
\put(40,10){\droite(-1,1){40}}
\put(40,10){\droite(1,1){40}}
\put(50,20){\droite(-1,1){30}}
\put(30,40){\droite(1,1){10}}
\put(70,40){\droite(-1,1){10}}
\endpicture}\kern 0.4ex}

\def\courbe{\kern-0.4ex
\hbox{\unitlength=.25pt
\picture(1600,420)(0,0)
\put(300,300){\droite(2,-1){100}}
\put(400,250){\droite(1,0){100}}
\put(500,250){\droite(2,1){100}}
\put(150,0){\droite(1,2){150}}
\put(0,400){\droite(3,-1){300}}
\put(250,320){$2=P$}

\put(0,200){\droite(0,1){200}}
\put(0,200){\droite(3,1){300}}
\put(0,200){\droite(3,-4){150}}
\put(150,0){\droite(1,0){300}}
\put(450,0){\droite(-1,2){150}}
\put(450,0){\droite(1,2){150}}
\put(300,300){\droite(2,1){100}}
\put(400,350){\droite(1,0){100}}
\put(500,350){\droite(2,-1){100}}

\put(0, 400){0}
\put(100, 0){0}
\put(460, 0){1}
\put(0, 150){1}
\put(600, 250){0}

\put(440, 150){$a$}
\put(440, 300){$b$}
\put(300, 100){$c$}
\put(140, 140){$d$}
\put(80, 300){$e$}

\put(900,300){\droite(2,-1){100}}
\put(1000,250){\droite(1,0){100}}
\put(1100,250){\droite(2,1){100}}
\put(750,0){\droite(1,2){150}}
\put(600,400){\droite(3,-1){300}}
\put(850,320){$1=P$}
\put(600, 400){0}
\put(700, 0){0}
\put(1200, 250){0}
\endpicture}\kern 0.4ex}

%op?ration gauche
\def\g{\mathop{\dashv}\nolimits}

%op?ration droite
\def\d{\mathop{\vdash}\nolimits}

%op?ration milieu
\def\m{\mathop{\perp}\nolimits}

%op?ration gauche dendriform
\def\l{\mathop{\prec}\nolimits}

%op?ration droite dendriform
\def\r{\mathop{\succ}\nolimits}

\def\cubeun{\kern-0.4ex
\hbox{\unitlength=1.5pt
\picture(200,80)(0,0)
\put(0,10){\droite(1,0){40}}
\put(32,30){\droite(1,0){38}}
\put(0,50){\droite(1,0){40}}
\put(0,48){\droite(1,0){40}}
\put(30,70){\droite(1,0){40}}

\put(0,10){\droite(0,1){38}}
\put(30,30){\droite(0,1){16}}
\put(32,30){\droite(0,1){16}}
\put(30,52){\droite(0,1){18}}
\put(32,52){\droite(0,1){18}}
\put(40,10){\droite(0,1){38}}
\put(70,30){\droite(0,1){38}}

\put(0,10){\droite(3,2){30}}
\put(2,10){\droite(3,2){30}}
\put(40,10){\droite(3,2){30}}
\put(42,10){\droite(3,2){28}}
\put(0,50){\droite(3,2){30}}
\put(40,50){\droite(3,2){30}}
\put(40,48){\droite(3,2){30}}

\put(140,20){\droite(1,0){40}}
\put(140,20){\droite(1,2){20}}
\put(180,20){\droite(-1,2){20}}

\put(0,5){$(\g)\g$}
\put(40,5){$(\d)\g$}
\put(0,43){$(\g)\d$}
\put(40,43){$(\d)\d$}
\put(70,63){$\d(\d)$}
\put(32,63){$\g(\d)$}
\put(27,25){$\g(\g)$}
\put(67,25){$\d(\g)$}

\put(100,30){$\mapsto$}

\put(130,10){$\g \g$}
\put(155,60){$\d \g$}
\put(180,10){$\d \d$}
\endpicture}\kern 0.4ex}

\def\cubedeux{\kern-0.4ex
\hbox{\unitlength=1.5pt
\picture(200,80)(0,0)
\put(0,10){\droite(1,0){40}}
\put(32,30){\droite(1,0){38}}
\put(0,50){\droite(1,0){40}}
\put(0,48){\droite(1,0){40}}
\put(30,70){\droite(1,0){40}}

\put(0,10){\droite(0,1){38}}
\put(30,30){\droite(0,1){16}}
\put(32,30){\droite(0,1){16}}
\put(30,52){\droite(0,1){18}}
\put(32,52){\droite(0,1){18}}
\put(40,10){\droite(0,1){38}}
\put(70,30){\droite(0,1){38}}

\put(0,10){\droite(3,2){30}}
\put(2,10){\droite(3,2){30}}
\put(40,10){\droite(3,2){30}}
\put(42,10){\droite(3,2){28}}
\put(0,50){\droite(3,2){30}}
\put(40,50){\droite(3,2){30}}
\put(40,48){\droite(3,2){30}}

\put(140,20){\droite(1,0){40}}
\put(140,20){\droite(1,2){20}}
\put(180,20){\droite(-1,2){20}}

\put(12,2){$(\m)\g$}
\put(41,35){$(\d)\m$}
\put(15,42){$(\m)\d$}
\put(1,30){$(\g)\m$}

\put(45,24){$\m(\g)$}
\put(31,54){$\g(\m)$}
\put(45,72){$\m(\d)$}
\put(72,45){$\d(\m)$}

\put(100,30){$\mapsto$}

\put(140,40){$\m \g$}
\put(170,40){$\d \m$}
\put(140,10){$(\g)\m = \m( \d)$}
%\put(154,30){$\m \m$}

\endpicture}\kern 0.4ex}

%fl\`eche projection
\def\proj{\mathop{\rightarrow \!\!\!\!\!\!\! \rightarrow}}

\def\xc{\check x}
 \def\x#1#2{#1^{(#2)}}

%\end

%This is trialgebra.tex in Plain TeX

\centerline {\grand Trialgebras and families of polytopes}
\B
\hfill

\hfill {\tt May 6, 2002}
\B

\centerline {\bf Jean-Louis Loday and Mar\' \i a O. Ronco} 
\BB
\N {\bf Abstract.} We show that the family of standard simplices and the family of Stasheff polytopes are dual to each other
 in the following sense.
 The chain modules of the standard simplices, resp. the Stasheff polytopes, assemble to give an operad. We show that these
 operads are
 dual of each other in the operadic sense. The main result of this paper is to show 
that they are both Koszul operads. As a consequence the generating series of the standard simplices and 
the generating series of the Stasheff polytopes are inverse to each other. The two operads give rise to new types of
algebras with 3 generating operations, 11 relations, respectively 7 relations, that we call {\it associative trialgebras} and 
{\it dendriform trialgebras}
respectively.  The free dendriform trialgebra, which is based on planar trees, has an interesting Hopf algebra structure, which
 will be dealt with in another paper. 

Similarly the family of cubes gives rise to an operad which happens to be self-dual for Koszul duality.
\BB

\N {\bf Introduction.} 
We introduce a new type of associative algebras characterized by the fact that the associative
product
$*$ is the sum of three binary operations :
$$x*y := x\l y + x\r y + x \cdot y \ ,$$
and that the associativity property of $*$ is a consequence of 7 relations satisfied by $\l , \r$ and
$\cdot$, cf. 2.1. Such an algebra is called a {\it dendriform trialgebra}. An example of a dendriform trialgebra is given by the algebra of
quasi-symmetric functions (cf. 2.3).

Our first result is to show that the free dendriform trialgebra on one generator can be described as an algebra over the
set of {\it planar trees}. Equivalently one can think of these linear generators as being the cells of the Stasheff
polytopes (associahedra), since there is a bijection between the $k$-cells of the Stasheff polytope of dimension $n$
and the planar trees with
$n+2$ leaves and $n-k$ internal vertices.

The knowledge of the free dendriform trialgebra permits us to construct the algebras over the dual operad 
(in the sense of 
 Ginzburg and Kapranov [G-K]) and therefore to construct the chain complex of a dendriform trialgebra. This dual type is
called the {\it associative trialgebra} since there is again three generating operations, and since all the relations are
of the associativity type (cf. 1.2). We show that the free associative trialgebra on one generator is linearly generated
by the cells of the standard simplices.

The main result of this paper is to show that the operads of dendriform trialgebras  (resp. associative trialgebras) is
a Koszul operad, or, equivalently, that the homology of the free dendriform trialgebra is trivial.

As a consequence of the description of the free trialgebras in the dendriform and associative framework, the
generating series of the associated operads are the generating series of the family of the Stasheff polytopes and of the
standard simplices respectively:
$$f^{\K }_t(x) = \sum_{n\geq 1} (-1)^n p(\K^{n-1},t)
x^n, \qquad f^{\DD }_t(x) = \sum_{n\geq 1} (-1)^n p(\DD^{n-1},t) x^n\ .$$
Here $p(X,t)$ denotes the Poincar\' e polynomial of the polytope $X$.

The acyclicity of the Koszul complex for the dendriform trialgebra operad implies that 
$$f^{\DD}_t(f^{\K}_t(x))=x.$$
Since $p(\DD^{n},t)=   ((1+t)^{n+1} -1)/t$ one gets 
$$ f^{\DD}_t(x) = {-x \over (1+x)(1+(1+t)x)}$$
and therefore
$$f^{\K}_t(x) = {-(1+(2+t)x) +\sqrt { 1+2(2+t)x+t^2x^2}\over 2(1+t)x}\ .$$

In [L1, L2] we dealt with dialgebras, that is with algebras defined by two generating operations. In the associative 
framework the dialgebra case is a quotient of the trialgebra case and in the dendriform framework the dialgebra case is 
a subcase  of the trialgebra case.
\M

If we split the associative relation for the operation  $*$ into 9 relations instead of 7, then we can devise a similar theory in which 
the family of Stasheff polytopes is replaced by the family of cubes. So we get a new type of algebras that we call the
{\it cubical trialgebras}. It turns out that the associated operad is self-dual (so the family of standard
simplices is to be replaced by the family of cubes). The generating series of this operad is  the  generating series of 
 the family of cubes:
$f^{I}_t(x) = {-x \over 1+(t+2)x}\ .$ It is immediate to check that $ f^{I}_t(f^{I}_t(x) ) = x\ ,$ hence one can presume
that this is a Koszul operad. Indeed we can prove that the Koszul complex of the cubical trialgebra operad is acyclic. 

As in the dialgebra case the associative algebra on planar trees can be endowed with a comultiplication which makes
it into a Hopf algebra. This comultiplication satisfies some compatibility properties with respect to the three
operations $\l, \r$ and $\cdot\  $. This subject will be dealt with in another paper.
\M

Here is the content of the paper.
\M

1. Associative trialgebras and standard simplices

2. Dendriform trialgebras and Stasheff polytopes

3. Homology and Koszul duality

4. Acyclicity of the Koszul complex

5. Cubical trialgebras and hypercubes

\M

In the first section we introduce the notion of associative trialgebra  and we compute the free algebra. This result gives the relationship with 
the family of standard simplices.

In the second section we introduce the notion of dendriform trialgebra  and we compute the free algebra, which is based on planar trees.
This result gives the relationship with  the  family of Stasheff polytopes.

In the third section we show that the associated operads are dual to each other for Koszul duality. Then we construct the chain complexes
which  compute the homology of these algebras. The acyclicity of the Koszul complex of the operad is equivalent to the acyclicity
of the chain complex of the free associative trialgebra. 

This acyclicity property is the main result of this paper, it is proved in the
fourth section. After a few manipulations involving the join of simplicial sets we reduce this theorem to proving the contractibility
 of some explicit simplicial complexes.
This is done  by producing a sequence of retractions by deformation. 

In the fifth section we treat the case of the family of hypercubes, along the same lines.

These results have been announced in [LR2]. 
\B

\BB

\N {\bf Convention.}  The category of vector spaces over the field $K$ is denoted by $\Vect$, 
and the tensor product of vector spaces over $K$ is denoted by $\t$. The symmetric group acting on $n$ elements is denoted by $S_n$. 

\BB

\N {\bf 1. Associative trialgebras and standard simplices.} 
\M

 In [L1, L2] the first author introduced the notion of associative dialgebra 
as follows. 
\M
\N {\bf 1.1 Definition.} An {\it associative dialgebra} is a vector space $A$ equipped with 2 binary operations : $\g$ called {\it left} 
and $\d$ called {\it right}, $$\leqalignno{
& \dashv\ : A\t A\to A,&{(\rm left)}\cr
& \vdash\ : A\t A\to A,&{(\rm right)}\cr} $$

satisfying the relations:
$$\left\{ \eqalign{ (x\g y)\g z &= x\g (y\g z) , \cr
 (x\g y)\g z &= x\g (y\d z) , \cr
 (x\d y)\g z &= x\d (y\g z) , \cr
 (x\g y)\d z &= x\d (y\d z) , \cr
 (x\d y)\d z &= x\d (y\d z) . \cr
 }\right.$$ 
Observe that the eight possible products with 3 variables $x,y,z$ 
(appearing in this order) occur in the relations. Identifying each product with a vertex of the cube and moding out the cube according 
to the relations transforms the cube into the triangle $\DD^2$ : 
$$\cubeun$$
The double lines indicate the vertices which are identified under the relations.

Let us now introduce a third 
operation $\m : A\t A \to A$ called {\it middle}. We think of left and right as being associated to the 0-cells of the interval and 
 middle to the 1-cell : 
$$\displaylines{ \g \qquad \m
\qquad \d \cr \bullet \!\! -\!\! -\!\! -\!\! -\!\! -\!\! -\!\! -\!\! -\!\! \bullet\cr }$$
Let us associate to any product in three variables a cell of the cube by using the three operations $\g , \d, \m $. The equivalence
 relation which transforms the cube into the triangle 
determines new relations (we indicate only the 1-cells):
 $$\cubedeux$$
This analysis justifies the following :
\M

\N {\bf 1.2 Definitions.} An {\it associative trialgebra} (resp. an {\it associative trioid}) is a vector space $A$ 
(resp. a set $X$) equipped with 3 binary operations : $\g$ called {\it left}, $\d$ called {\it right} and $\m$ called 
{\it middle}, satisfying the following 11 relations :
$$\left\{ \eqalign{
(x\g y)\g z &= x\g (y\g z) , \cr
(x\g y)\g z &= x\g (y\d z) , \cr
(x\d y)\g z &= x\d (y\g z) , \cr
(x\g y)\d z &= x\d (y\d z) , \cr
(x\d y)\d z &= x\d (y\d z) , \cr
}\right.
$$
$$
\left\{ 
\eqalign{
(x\g y)\g z &= x\g (y\m z) , \cr
(x\m y)\g z &= x\m (y\g z) , \cr
(x\g y)\m z &= x\m (y\d z) , \cr
(x\d y)\m z &= x\d (y\m z) , \cr
(x\m y)\d z &= x\d (y\d z) , \cr
}\right.
$$
$$
\left\{ \eqalign{
(x\m y)\m z &= x\m (y\m z) . \cr
}\right.$$

First, observe
that each operation is associative. Second, observe that the following rule holds: ``on the bar side, does not matter 
which product''. Third, each relation has its symmetric counterpart which consists in reversing the order of the 
parenthesizing, exchanging $\d$ and $\g$, leaving $\m$ unchanged. 

A morphism between two associative trialgebras is a linear map which is compatible with the three operations. 
We denote by {\bf Trias} the category of associative trialgebras. 
\M

\N {\bf 1.3 Relationship
with the planar trees.} The set of planar trees with $(n+1)$ leaves is denoted by $T_n$, see 2.4 for notation 
and definitions.
We associate the
trees in $T_2$ to the three binary operations as follows: $$\displaylines{
\big(\ \arbreC ; x,y\big) \mapsto \big(\ \arbreA ; x\g y\big)\cr 
\big(\ \arbreB ; x,y\big) \mapsto \big(\ \arbreA ; x\d y\big)\cr 
\big(\ \arbreBC ; x,y\big) \mapsto \big(\ \arbreA ; x\m y\big).\cr }$$
Observe that it is the direction of the middle leaf which determines the operation. Anyone of the 11 trees $t$ in
$T_3$ gives two different ways of computing the image of $(t; x,y,z)$. 
Equating the two results gives a relation. For instance, let $t=\arbreuut $. The first computation gives $$
(\ \arbreuut ; x,y,z) \mapsto (\ \arbreB ; x\m y,z) \mapsto (\ \arbreA; (x\m y )\d z\ ). $$
The second computation gives
$$
(\ \arbreuut ; x,y,z) \mapsto (\ \arbreB ; x,y\d z) \mapsto (\ \arbreA; x\d (y \d z)). $$
So this tree gives rise to the 10th relation of the list 1.2. It is straightforward to verify that the 
11 trees of $T_3$ give the 11 relations of 1.2. This relationship will be exploited in constructing the chain complex of an associative
trialgebra in section 3.
\M

\N {\bf 1.4 Examples of associative trialgebras.} 
\S
\N (a) If $A$ is an associative trialgebra, then the $n\times n$-matrices over $A$ still form an associative trialgebra
 by taking the operations coefficient-wise. 

\N (b) If $\g = \m = \d$, then we get simply an associative algebra (nonunital). So we get a functor between the categories 
of algebras : $$ {\bf As} \to {\bf Trias}\ . $$ 
Ignoring the operation $\m$ gives an associative dialgebra. Hence there is a (forgetful) functor 
$${\bf Trias} \to {\bf Dias}$$
 from the category of trialgebras to the category of dialgebras. 
\S
\N (c) The vector space over an associative trioid is obviously an associative trialgebra.
\S
\N d) {\it The Solomon algebra.} Let $V=\oplus _{n\geq 0}K\cdot \omega _n$ be the graded $K$-vector space such that the
subspace of homogeneous elements of degree $n$ is the vector space of dimension one, spanned by the generator $\omega _n$, 
for all $n\geq 0$. Consider the tensor algebra $T(V)$, with the operations $\perp $, $\dashv $ and $\vdash $ given by:
$$(\omega _{n_1}\otimes \dots \otimes \omega _{n_r})\perp (\omega _{m_1}\otimes \dots \otimes \omega _{m_k}):=
\omega _{n_1}\otimes \dots \otimes \omega _{n_r}\otimes \omega _{m_1}\otimes \dots \otimes \omega _{m_k},$$
$$(\omega _{n_1}\otimes \dots \otimes \omega _{n_r})\dashv (\omega _{m_1}\otimes \dots \otimes \omega _{m_k}):=
\omega _{n_1}\otimes \dots \otimes \omega _{n_r}\otimes \omega _{m_1+\dots +m_k},$$
$$(\omega _{n_1}\otimes \dots \otimes \omega _{n_r})\vdash (\omega _{m_1}\otimes \dots \otimes \omega _{m_k}):=
\omega _{n_1+\dots +n_r}\otimes \omega _{m_1}\otimes \dots \otimes \omega _{m_k},$$
for $n_1,\dots ,n_r,m_1,\dots ,m_k\geq 0$.
It is easy to check that $(T(V), \perp ,\dashv ,\vdash )$ is an associative trialgebra. The associative algebra $(T(V),\perp
)$ is isomorphic to the Solomon algebra $Sol_{\infty }$ (cf. for instance [LR1]).

\M

\N {\bf 1.5 Notation}. Let $[n-1]:= \{0, \cdots , n-1\}$ be a set with $n$ elements. The set of non-empty subsets of $[n-1]$ 
is denoted by $P_n$. Observe that $P_n$ is graded by the cardinality of its members. We denote by $P_{n,k}$ the subset of $P_n$ 
whose members have cardinality $k$. So $P_n= P_{n,1}\cup \cdots \cup P_{n,n}$. \M

\N {\bf 1.6 Free associative trialgebra.} By definition the free associative trialgebra  over the vector space $V$ is an
 associative trialgebra  $Trias(V)$ equipped
with a map $V\to Trias(V)$, which satisfies the following universal property. For any map $V\to A$, where $A$
 is an associative trialgebra, 
there is a unique extension $Trias(V) \to A$ which is a morphism of associative trialgebras. 

Since the operations have no symmetry and since the relations let the variables in the same order, $Trias(V)$ is 
completely determined by the free associative trialgebra on one generator (i.e. $V=K$). The latter is a
graded vector space of the form
$$ Trias (K) = \oplus_{n\geq 1} Trias (n). $$ 
From our motivation of
defining the associative trialgebra type it is clear that for $n=1, 2, 3$, a basis of $Trias (n)$ is given by the 
elements of $P_1$, $P_2$ and $P_3$ respectively (i.e. the cells of $\DD ^0, \DD^1, \DD^2$ respectively). 
\M

Let us denote by
$$ {\rm bij} : [{i_1-1}]\cup \cdots \cup [ {i_n-1}] \to [{i_1+ \cdots +i_n-1}] $$ 
the bijection which sends $k\in [i_j -1]$ to $i_1+ \cdots +i_{j-1}+k \in [i_1+ \cdots +i_n-1]$. 
\M

\N {\bf 1.7 Theorem.} {\it The free associative trialgebra $Trias (K)$ on one generator is  $\oplus_{n\geq 1} K[P_n]$ as a vector space.
The binary operations $\g , \m $ and $\d$ from 
$K[P_p]\t K[P_q]$ to $K[P_{p+q}]$ are given by 
$$ X\g Y = {\rm bij} (X) , \quad X\m Y =
 {\rm bij}(X \cup Y) \quad X\d Y= {\rm bij} (Y)\ , 
$$
 where $X\in P_p$ and $Y\in P_q$ and  ${\rm bij} : [p-1]\times [q-1] \to [p+q-1] $. } 
\M

\N {\bf 1.8 Corollary.} {\it The free associative trialgebra $Trias (V)$ on the vector space $V$ is 
$$ Trias (V) = \oplus _{n\geq 1}K[P_n]\t V^{\t n},$$ 
and the operations are induced by the operations on $Trias (K)$ and concatenation.}
\M
\N {\it Proof.} It suffices to make explicit the free trioid in one generator, see Proposition 1.9 below. Indeed, it proves Theorem 1.7 
by applying the functor which sends a set $Z$ to the vector space $K[Z]$ having the elements of $Z$ as a basis. Then the
Corollary is a consequence of the Theorem because all the relations in the definition of an associative trialgebra leave
the variables in the same order. \hfill $\square$
\M

\N {\bf 1.9 Proposition.} {\it The free trioid ${\cal T}$ on one generator $x$ is isomorphic to the trioid $P =
\bigcup_{n\geq 1} P_n$ equipped with the operations described in Theorem 1.7 above.}
\M
\N {\it Proof.} First we  prove that $(P ; \g , \d , \m )$ is a trioid generated by $\{ 0\} \in
P_1$. For convenience let us denote this generator by $x$ and by $x\cdots \xc \cdots \xc \cdots x$ the  element corresponding to $X\in
P_n$, where there are
$n$ copies of $x$ and, if $i\in X$, then the $i$th factor  is ceched. For instance $\{0,2\}\in P_5$ corresponds to $\xc x\xc xx$. 
Under this notation the operations are easy to describe: one concatenates
the two elements, keeping only the marking  on the left side for $\g$, on the right side for $\d$, on both sides for $\m$. For instance
$$\left\{ \eqalign{
\xc x\g x\xc x = \xc xxxx,\cr
\xc x\d x\xc x = xxx\xc x,\cr
\xc x\m x\xc x = \xc xx\xc x.\cr
}\right.$$
It is immediate to verify that the eleven relations are fulfilled.

Since ${\cal T}$ is the free trioid generated by $x$, there exists a unique trioid morphism $\phi : {\cal T} \to P$. Each map  $\phi_n : {\cal
T}_n \to P_n$ is surjective since, in $P$, $x=\{0\}\in P_1$ is also the generator. In order to prove that $\phi $ is an isomorphism,
 it suffices to show
that $\# {\cal T}_n \leq \# P_n$.
\M

\N{\bf 1.10 Lemma.} {\it Any complete parenthesizing of
$${\underbrace{ (x\d \cdots \d x) } \atop a_0}
{ \d  \atop }
 {\underbrace{ (x\g \cdots \g x) }\atop a_1}
{ \m  \atop }
 {\underbrace{ (x\g \cdots \g x) }\atop a_2}
{ \m \cdots \m  \atop }
 {\underbrace{ (x\g \cdots \g x)  }\atop a_k}
$$
where $a_0 \geq 0, a_i \geq 1$ for $i= 1, \cdots k$, gives the same element, denoted $\oo$, in ${\cal T}$. We call it the normal form of
$\oo$. Its image under $\phi$ in $P$ is 
$${\underbrace{ x\cdots x} \atop a_0} {\underbrace{\xc\cdots x }\atop a_1}{\underbrace{\xc\cdots x }\atop a_2}{ \cdots x \atop }
{\underbrace{\xc\cdots x } \atop a_k}.$$}

\N {\it Proof.} Putting parentheses outside (resp. inside) the existing parentheses does not change the value of the element by virtue of
relations 9 and 11 (resp. 1 and 5). The second statement is immediate by direct inspection. \hfill $\square$
\M
\N {\it End of the proof of Proposition 1.9.} 
Since any element in $P_n$ is the image of an element of the type indicated in Lemma 1.10, it suffices to show that any element
 in ${\cal T}_n $ can be written under this form. We work by induction on $n$. It is clear for $n=1$. We suppose that it is true for all $p< n$.
Any $\ss \in {\cal T}_n$ is of the form $\ss ' \g \ss '' $ or $\ss ' \d \ss '' $ or $\ss ' \m \ss '' $ for some $\ss ' \in {\cal T}_p, 
\ss ''\in {\cal T}_q $. We write $\ss '$ and $\ss ''$ in a normal form as in Lemma 1.10 and we compute the three elements  $\ss ' \g \ss ''
,\ss ' \d \ss '' $ and $\ss ' \m \ss '' $. By using the relations 1 to 11 it is easy to show that they can be written under a normal
 form. So the proof of Proposition 1.9 is complete.
\hfill $\square$
\M

\N {\bf 1.11 Filtration.} The set $P_n$ can be filtered by $F_kP_n := \cup_{i\leq k} P_{n,i}$, cf. 1.5. Since, in any product of two elements, the
number of  marked variables is equal or less than the sum of the numbers of the components, the image of $F_kP_n\times F_lP_m$ is in
$F_{k+l}P_{n+m}$.
\M

\N {\bf 1.12 The family of standard simplices.} Let $\DD^n = \{ (x_0, \cdots , x_n) \in {\bf R}^{n+1} \mid x_0+\cdots + x_n= 1, 0\le x_i
\le 1\} $
 be the standard $n$-simplex. As usual we label its vertices by the integers $0$ to $n$. So the vertex $
 i$ has coordinates 0 except $x_i = 1$. An $i$-cell in $\DD^n$ is completely determined by its vertices, hence by a non-empty subset of 
$[n] = \{0, \cdots , n\}$, that is, following our notation, an element of $P_{n+1}$. So there is a bijection between the 
$k$-cells of $\Delta ^n$ and the set $P _{n+1,k+1}$. 
$$\simplexes$$
Observe that the Poincar\'e polynomial of $\DD^n$ is 
$$p(\DD^n, t) :=
\sum_{k\geq 0} \# (k\hbox {-cells}) t^k =  {((1+t)^{n+1} -1)\over t}.$$
If we define the generating series of a family of polytopes $X(n), n\geq 0$ by 
$f^{X }_t(x) = \sum_{n\geq 1} (-1)^n p(X(n-1),t) x^n$, then we get the following for the family of standard simplices: 
$$ f^{\DD}_t(x) = {-x \over (1+x)(1+(1+t)x)}\ .$$ 
\M

\N {\bf 1.13 Generating series of a filtered operad.} The operad $\P$ determined by a category of algebras is a functor
$\P : \Vect \to \Vect$ of the form $\P(V) = \oplus_{n\geq 1} \P(n)\t_{S_n} \tT Vn $ (here $\P(n)$ is a right
$S_n$-module) together with an  associative and unital transformation of functors $\cc : \P \circ \P \to \P$, cf. [G-K, O]. The free
$\P$-algebra over $V$  is precisely $\P(V)$. 

By definition the generating series of an operad $\P$ is 
$$f^{\P }(x) := \sum_{n\geq 1} (-1)^n (1/{n!}) \dim \P (n) x^n.
$$ 
When the operad is filtered, we can define a finer invariant by replacing the dimension of $\P (n)$ by its
Poincar\'e  polynomial 
$$p(\P (n),t):=\sum_{k\geq 0} \dim (F_k\P _n/F_{k-1}\P _n)\ t^k$$
to get a series with polynomial in $t$ as coefficients: 
$$f^{\P }_t(x) := \sum_{n\geq 1} (-1)^n (1/{n!})p(\P(n),t) x^n.$$
In [G-K] it is shown that if the  quadratic operad $\P$ is a Koszul operad, then the generating series of $\P$ and of its dual
  $\P^!$ are related by $ f^{\P}(f^{\P ^!}(x) ) = x\ $. This formula is obtained by computing the Euler-Poincar\'e characteristic
 of the Koszul complex of $\P$, which gives the left hand-side. Since this complex is acyclic, its homology is trivial, and this
 gives the right-hand side. 

If the quadratic operad $\P$ is filtered, then a refinement of this argument gives the functional equation: 
$$ f^{\P}_t(f^{\P ^!}_t(x) ) = x\ .$$
\M

\N {\bf 1.14 The operad of associative trialgebras.} 
Let $Trias$ be the operad associated to the associative trialgebras. By Corollary 1.8 we have $Trias(n) = K[P_n]\t K[S_n]$, 
where $S_n$ is the symmetric group. So, as an $S_n$-module, $Trias(n) $ is the direct sum of several copies of the regular representation, one 
for each element in $P_n$. In other words $Trias$ is a non-$\SS$-operad in the sense of [O, p.4]. The filtration of the free trioid 
described in 1.11 induces a filtration on
  the operad  $Trias$. Since this filtration corresponds precisely to the filtration of the standard simplex by the dimension of
 the cells, the generating series are equal:
 
$$ f^{Trias }_t(x) = f^{\DD}_t(x) = {-x \over (1+x)(1+(1+t)x)}\ .$$ 

\N {\bf 1.15 Relationship with the Leibniz and Poisson algebra structures.} The notion of associative dialgebra was first 
introduced as an analogue of associative algebra for Leibniz algebras. Let us recall that a Leibniz algebra is defined by a
 binary operation $[-,-]$ which is not necessarily skew-symmetric and satisfies the right Leibniz identity:
$$ [[x,y],z] = [[x,z],y] + [x,[y,z]].$$
If the bracket happens to be skew-symmetric, then this is a Lie bracket. Any associative dialgebra gives rise to a Leibniz bracket by:
$$ [x,y]:= x\g y - y\d x.$$
Suppose now that we would like to construct a {\it noncommutative} version of {\it Poisson algebra}. Then we introduce an associative
 operation $xy$ (not necessarily commutative), and it is natural to require that its relationship with the Leibniz bracket is given by
$$\eqalignno{
[xy,z] &= x[y,z] + [x,z] y\ ,& (1.15.1)\cr
[x, yz-zy] &= [x, [y,z]]\ .& (1.15.2)\cr
}$$

\N {\bf 1.16 Proposition.} {\it Let $(A, \g, \d,\m)$ be an associative trialgebra. By defining
$$[x,y]:= x\g y - y\d x \quad {\rm and }\quad xy:= x\m y $$
we get a noncommutative Poisson algebra structure on $A$.}
\M

\N {\it Proof.} The relation 1.15.1 is a consequence of the relations number 7, 8 and 9 in 1.2 and the relation 1.15.2 is a consequence of
 the relations number 6 and 10.\hfill $\square$

Compare with the work of Marcelo Aguiar [A].
\M
\N {\bf 1.17 Relationship with the boundary map of the standard simplex.} The space $K[P_n]$ is in fact the chain module of the standard simplex
$\DD^{n-1}$, and so it is equipped with a differential map $\dd : K[P_{n,k}]\to K[P_{n, k-1}]$. Explicitly $\dd$ is given by 
$$\delta (X):=\sum _{i=1}^r(-1)^{i+1}X\setminus \lbrace n_i\rbrace ,$$
for $X=\lbrace n_1<n_2<\dots <n_k\rbrace $ a subset of $[n-1]$. 
The relationship of $\dd$ with the three operations $\g, \d$ and $\m$ is given (for\ $X\in P_{n,k}$) by:

$$\eqalign{
\delta (X\dashv Y)&=\delta (X)\dashv Y, \cr
\delta (X\vdash Y)&=(-1)^k X\vdash \delta (Y),\cr
\delta (X\perp Y)&=\cases {
\delta (X)\perp Y + (-1)^kX\perp \delta (Y)& \ for\ $\delta (X)\neq 0$ and $\delta (Y)\neq 0$,\cr
X\vdash Y+(-1)^kX\perp \delta (Y)&\ for\  $\delta (X)=0$ and $\delta (Y)\neq 0$,\cr
\delta (X)\perp Y+(-1)^kX\dashv Y&\ for\  $\delta (X)\neq 0$ and $\delta (Y)=0$,\cr
X\vdash Y - X\dashv Y&\ for\  $\delta (X)=0$ and $\delta (Y)= 0$.\cr }
\cr}$$
\M

\N {\bf 2. Dendriform trialgebras and Stasheff polytopes.} In [L1, L2] the first author introduced the notion of
dendriform dialgebras.  Here we add a third operation. 
\M

\N {\bf 2.1 Dendriform trialgebras.} By definition a {\it dendriform trialgebra  } is a vector space $D$ equipped with three binary operations :
\S
$\l$ called {\it left},\quad   $\r$ called {\it right},\quad $\cdot $ called {\it middle}, \S
satisfying the following relations :
$$\left\{ \eqalign {
(x \l y) \l z &= x \l (y * z)\ , \cr
(x \r y) \l z &= x \r (y \l z)\ , \cr
(x * y) \r z &= x \r (y \r z)\ , \cr
}\right.
$$
$$
\left\{ \eqalign {
(x \r y) \cdot z &= x \r (y \cdot z)\ , \cr (x \l y) \cdot z &= x \cdot (y \r z)\ , \cr (x \cdot y) \l z &= x \cdot (y \l z)\ , \cr 
}\right.
$$
$$
\left\{ \eqalign {
(x \cdot y) \cdot z &= x \cdot (y \cdot z)\ , \cr }\right.
$$
where $x*y := x\l y + x\r y  + x\cdot y$. 
\M

\N {\bf 2.2 Lemma. } {\it The operation $*$ is associative.} 
\M

\N {\it Proof.} It suffices to add up all the relations to observe that on the right side we get $(x*y)*z$ and on the left side $x*(y*z)$. 
Whence the assertion. \hfill $\square$ 
\M

In other words, a dendriform trialgebra is an associative algebra for which the associative operation is the sum of 
three operations and the associative relation splits into 7 relations.

 We denote by ${\bf Tridend}$ the category of
dendriform trialgebras and by
$Tridend$ the associated operad. By the preceding lemma, there is a well-defined functorÊ:
$$ {\bf Tridend} \to {\bf As}\ ,$$
where {\bf As} is the category of (nonunital) associative algebras. 

Observe that the operad $Tridend $ does not come from a set operad because the operation $*$ needs a sum to be defined. However there
is a property which is close to it. It is discussed and exploited in [L3].
\M

\N {\bf 2.3 Examples of dendriform trialgebras.}
\M
\N (a) If $D$ is a dendriform trialgebra, then the
$n\times n$-matrices over $D$ still form a dendriform trialgebra. 
\S
\N (b) If the operation $\cdot$ is taken to be trivial (i.e $x\cdot y =0$ for any $x,y\in D$), then $x * y = x\l y + x\r
y$, and we get simply a dendriform dialgebra as defined in [L1, L2] (two generating operations and 3 relations).
This defines a functor 
$${\bf Didend}\to {\bf Tridend}.$$

\N (c) {\it Quasi-symmetric functions.} Let $K<y_1, y_2, \cdots , y_k, \cdots >$ be the free associative unital algebra on a countable set of
variables
$y_k$. We define a new associative product on it by the following inductive formula
$$y_k \oo * y_{k'} \oo ' := y_k (\oo * y_{k'} \oo ') +  y_{k'} (y_k \oo *  \oo ') + y_{k+k'}( \oo *  \oo ')\ ,$$ 
where $\oo$ and $\oo '$ are monomials or $1$ (unit for $*$). So for instance
$$ y_k  * y_{k'}  := y_k  y_{k'}  +  y_{k'} y_k  + y_{k+k'} \ .$$ 
If we denote by $y_k \oo\l y_{k'} \oo '$, resp.  $y_k \oo \r y_{k'} \oo '$, resp.  $y_k\oo \cdot y_{k'} \oo '$, the
first, resp. second, resp. third summand in this sum, then we can show that we have defined a dendriform trialgebra
structure on the augmentation ideal. Indeed, let $x=y_k \oo\ , y= y_l \oo '$ and $z= y_m \oo ''$. Then the seven relations of 2.1 hold
and give the elements
$$\displaylines{
y_k ( \oo * y_l \oo ' * y_m \oo '')\ , \cr
y_l ( y_k \oo * \oo ' * y_m \oo '')\ , \cr
y_m( y_k \oo * y_l \oo ' *  \oo '')\ , \cr
y_{k+l} ( \oo *  \oo ' * y_m \oo '')\ , \cr
y_{l+m} ( y_k \oo *   \oo ' *   \oo '')\ , \cr
y_{k+m} ( \oo * y_l \oo ' *  \oo '')\ , \cr
y_{k+l+m} ( \oo *  \oo ' *   \oo '')\ . \cr
}$$
Equipped with the associative (and commutative) product $*$,  the space of noncommutative polynomials $K\langle 
y_1, y_2, \cdots , y_k, \cdots \rangle $ is the algebra ${\rm QSym}$ of
 quasi-symmetric functions, cf. [C, formula 94], [H].

One can show that the augmentation ideal is the free dendriform commutative trialgebra on one generator $y_1$. Here commutative
 means $x\r y = y\l x$ and $x\cdot y = y\cdot x$ for any $x$ and $y$.
\S

(d) {\it The partition example.} For $n\geq 1$, a partition of
$n$ is a family of positive integers
${\underline {\bf n}}=(n_1,\dots ,n_r)$, such that $\sum _{i=1}^rn_i = n$. Given partitions ${\underline {\bf n}}= 
(n_1,\dots ,n_r)$ of $n$ and ${\underline {\bf m}}= (m_1,\dots ,m_k)$ of $m$, we denote by ${\underline {\bf n}}\times 
{\underline {\bf m}}$ the partitions of $n+m$ given by:
$${\underline {\bf n}}\times {\underline {\bf m}}:=(n_1,\dots ,n_r,m_1,\dots ,m_k).$$
Given a partition ${\underline {\bf n}}=(n_1,\dots ,n_r)$ of $n$, a {\it $(n_1,\dots ,n_r)$-shuffle} is a permutation $\sigma $
in the symmetric group $S_n$, such that:
$$\displaylines{
\sigma (1)<\dots <\sigma (n_1),\ \sigma (n_1+1)<\dots <\sigma (n_1+n_2),\ \dots ,\hfill \cr
\hfill  \sigma (n_1+\dots +n_{r-1}+1)<\dots <\sigma (n).\cr
}$$
We denote by $Sh(n_1,\dots ,n_r)$ the set of all $(n_1,\dots ,n_r)$-shuffles.

\N Given permutations $\sigma \in S_n$ and $\tau \in S_m$, we denote by $\sigma \times \tau$ the element of $S_{n+m}$ whose
image is $(\sigma (1),\dots ,\sigma (n),\tau (1)+n,\dots ,\tau (m)+n)$.

\N Given an $(n_1,\dots ,n_r)$-shuffle $\sigma $ there exist
unique elements $\sigma _0\in Sh(n_1+\dots +n_{r-1},n_r)$ and $\sigma _1\in Sh(n_1,\dots ,n_{r-1})$, such that:
$$\sigma = \sigma _0(\sigma _1\times 1_{S_{n_r}}),\eqno (*)$$
where $1_{S_{n_r}}$ is the identity element of the group $S_{n_r}$.
\M

For $n\geq 1$, the set of ordered partitions of $n$ is the set
$$\Pi _n:=\lbrace ({\underline {\bf n}},\sigma ):\ {\rm where}\ {\underline {\bf n}}\ {\rm is\ a\ partition\ of}\ n\ {\rm
and}\ \sigma \in Sh(n_1,\dots ,n_r)\rbrace .$$
If ${\underline {\bf n}} = (1, 2, \cdots , n)$, then we simply denote it by $(n)$. The element $((n), 1_{S_n})$ will play a particular role.

Given elements $({\underline {\bf n}},\sigma )\in \Pi _n$ and $({\underline {\bf m}},\tau )\in \Pi _m$, and an $(n,m)$-shuffle
$\gamma $, we define a new element $({\underline {\bf n}},\sigma )\times _{\gamma }({\underline {\bf m}},\tau )$ in $\Pi
_{n+m}$ as follows:
$$
({\underline {\bf n}},\sigma )\times _{\gamma }({\underline {\bf m}},\tau ):=({\underline {\bf n}}\times {\underline {\bf
m}}, \gamma (\sigma \times \tau )).
$$
\M

Let $K[\Pi _{\infty }]$ be the graded vector space spanned by the graded set $\cup _{n\geq 1}\Pi _n$. 
\N Given partitions ${\underline {\bf n}}= (n_1,\dots ,n_r)$ of $n$ and ${\underline {\bf m}}= (m_1,\dots ,m_k)$ of $m$, and
permutations $\sigma \in Sh(n_1,\dots ,n_r)$ and $\tau \in Sh(m_1,\dots ,m_k)$, we define the three operations on $K[\Pi _{\infty }]$
 as follows.
\S 
\N {\it $\bullet$ Right operation:}  If $({\underline {\bf m}},\sigma )=((m),1_{S_m})$, then:
$$({\underline {\bf n}},\sigma )\succ ((m),1_{S_m}):=((n_1,\dots ,n_r,m),\sigma \times 1_{S_m}).$$
\S

\N If $k\geq 2$, then
$$({\underline {\bf n}},\sigma )\succ ({\underline {\bf m}},\tau ):=(({\underline {\bf n}},\sigma )*({\underline {\bf
m\rq }},\tau _1))\times _{1_{S_n}\times \tau _0}((m_k),1_{S_{m_k}}),$$
where $\tau _0$ and $\tau _1$ are the permutations defined in formula  $(*)$ above and ${\underline {\bf m\rq }}=(m_1,\dots ,m_{k-1})$.
\M

\N {\it  $\bullet$ Left operation:} If $r=1$, then $((n),\sigma )=((n),1_{S_n})$. In this case, we define
$$((n),1_{S_n})\prec ({\underline {\bf m}},\tau ):=({\underline {\bf m}},\tau )\times _{\alpha _{n,m}}((n),1_{S_n})$$
 where $\alpha
_{n,m}(i):=\cases { n+i&\ for\ $1\leq i\leq m$\cr
i-m&\ for\ $m+1\leq i\leq m+n$.\cr }$
\S

\N If $r\geq 2$, then we have
$$(({\underline {\bf n}},\sigma )\prec ({\underline {\bf m}},\tau ):= (({\underline {\bf
n\rq }},\sigma _1 )*({\underline {\bf m}},\tau ))\times _{\beta }((n_r),1_{S_{n_r}}),$$
where 
$$\beta (i):=\cases {
\sigma _0(i)&\ for\ $1\leq i\leq n_1+\dots +n_{r-1}$,\cr
i+n_r&\ for\ $n_1+\dots +n_{r-1}<i\leq n_1+\dots +n_{r-1}+m$\cr
\sigma _0(i-n_1-\dots -n_{r-1})&\ for\ $n_1+\dots +n_{r-1}+m<i\leq n+m.$\cr },$$
\N and $\sigma _0$ and $\sigma _1$ are the elements defined in formula  $(*)$ above.
\M 

\N {\it  $\bullet$ Middle operation:} If $({\underline {\bf n}},\sigma )=((n),1_{S_n})$ and $({\underline {\bf m}},\tau )=((m),1_{S_m})$, then
$$((n),1_{S_n})\cdot ((m),1_{S_m}):=((n+m),1_{S_{n+m}}).$$
\S

\N If $({\underline {\bf n}},\sigma )=((n),1_{S_n})$ and $k\geq 2$, we have:
$$((n),1_{S_n})\cdot ({\underline {\bf m}},\tau ):=((m_1,\dots ,m_{k-1},m_k+n),\beta (\tau _1\times 1_{S_{m_k+n}})),$$
where 
$$\beta :=\cases {
\tau _0(i)&\ for\ $1\leq i\leq m_1+\dots +m_{k-1}$,\cr
i-m_1-\dots -m_{k-1}&\ for\ $m_1+\dots +m_{k-1}<i\leq m_1+\dots +m_{k-1}+n$,\cr
\tau _0(i-n)&\ for\ $m_1+\dots +m_{k-1}+n<i\leq m+n$.\cr }$$
\S

\N if $({\underline {\bf m}},\tau )=((m),1_{S_m})$, and $r\geq 2$, then:
$$({\underline {\bf n}},\sigma )\cdot ((m),1_{S_m}):=((n_1,\dots ,n_{r-1},n_r+m),\sigma \times 1_{S_m}).$$
\S

\N If $k,r\geq 2$, then the product $\cdot $ is given by:
$$({\underline {\bf n}},\sigma )\cdot ({\underline {\bf m}},\tau ):=(({\underline {\bf n\rq }},\sigma _1 )*
({\underline {\bf m\rq }},\tau _1))\times _{\alpha }((n_r+m_k), 1_{S_{n-r+m_k}}),$$
where 
$$\alpha (i):=\cases {
\sigma _0(i)&\ for\ $1\leq i\leq n_1+\dots +n_{r-1}$,\cr
\tau _0(i-n_1-\dots -n_{r-1})&\ for\ $n_1+\dots +n_{r-1}<i\leq n+m-n_r-m_k$,\cr
\sigma _0(i-m_i-\dots -m_{k-1})&\ for\ $n+m-n_r-m_k<i\leq n+m-m_k$,\cr
\tau _0(i-n)&\ for\ $n+m-m_k<i\leq n+m$.\cr }$$

With these definitions one can show that $(K[\Pi _{\infty }]; \prec, \succ, \cdot)$ is a dendriform trialgebra.
\M

\N {\bf 2.4 Planar trees.} We denote by $T_n$ the set of planar trees with $n+1$ leaves, $n\geq 0$ (and one root) such that the valence of each internal vertex
is at least 2. 
Here are the first of them:
$$T_0 = \{\vert\} ,\qquad \ T_1 = \{\ \arbreA\} ,\qquad T_2 = \{\ \arbreB , \arbreC ; \arbreBC\}$$ 
$$T_3 = \{\ \arbreun, \arbredeux, \arbretrois, \arbrequatre, \arbrecinq; \arbreuut, \arbretut, \arbretuu, \arbreutt, \arbrettu; \arbrettt \}. $$
The integer $n$ is called the {\it degree} of $t\in T_n$.
The number of elements in $T_n$ is the so-called {\it super Catalan number} $C_n$ :
\M

{\vbox{
\hrule
\halign{&\vrule#&\strut \hfil# \cr
$n$	         && 1 && 2 && 3		&&  4 && 5 && \cr
\noalign{\hrule}
$C_n$	&& 1 && 3 && 11		&& 45 && 197 && \cr
\noalign{\hrule}
  }}
\S

 The set $T_n$ is the disjoint union of the
sets $T_{n,k}$  made of the planar trees which have $n-k+1$ internal vertices. For instance $T_{n,1}$ is made of the
planar binary trees, and  its cardinality is the Catalan number $(2n)!/{n!(n+1)!}$. On the other extreme the set $T_{n,n}$ has only one
 element, which is the planar tree with one vertex. It is sometimes called a {\it corolla}.  So we have 
$$T_n = T_{n,1} \cup \cdots \cup T_{n,n}.$$
By convention $T_0 = T_{0,0}$.

The grafting of $k$ planar trees $x^{(0)}, \cdots , x^{(k)}$ is a planar tree denoted $x^{(0)}\vee \cdots
\vee x^{(k)}$ obtained by joining the $k+1$ roots to a new vertex and adding a new root. Any planar tree can be uniquely obtained as
$x=x^{(0)}\vee \cdots \vee x^{(k)}$, where $k+1$ is the valence of the lowest vertex.
We will use the uniqueness of this decomposition in the construction of a dendriform trialgebra structure on planar trees. Observe
 that the degree of $x^{(i)}$ is strictly smaller than the degree of $x$.
\M

\N {\bf 2.5 Free dendriform trialgebra.} The free dendriform trialgebra  over the vector space $V$ is a dendriform trialgebra 
$Tridend(V)$ equipped with a map $V\to Tridend(V)$  which satisfies the classical universal property, cf. 1.6. In the following theorem
we make it explicit in terms of planar trees. 
\M

\N {\bf 2.6 Theorem.} {\it The free dendriform trialgebra on one generator is 
$$Tridend (K) = \oplus_{n\geq 1} K[T_n],$$ 
where $T_n$ is the set of planar trees with $(n+1)$ leaves. 

The binary operations are given on
$T_p\times T_q$ by the recursive formulas: 
$$\eqalign{
x\l y &= x^{(0)}\vee \cdots \vee (x^{(k)}* y)\ , \cr x\cdot y &= x^{(0)}\vee \cdots \vee (x^{(k)}*y^{(0)})\vee \cdots \vee y^{(\ell )}\ ,
 \cr x\r y &= (x*y^{(0)})\vee \cdots \vee y^{(\ell )}\ , \cr }$$
where $x = x^{(0)}\vee \cdots \vee x^{(k )}\in T_p$ and $y = y^{(0)}\vee \cdots \vee y^{(\ell )}\in T_q.$
As before $x*y := x\l y + x\cdot y + x\r y$ and $\vert \in T_0$ is a unit for $*$.} 
\M

\N {\it Proof.} It follows from the following two lemmas. In the first one we prove that  $(\oplus_{n\geq 1} K[T_n]; \l , \r , \cdot )$  is a
dendriform trialgebra generated by the tree $\arbreA$. As a consequence there is a unique dendriform trialgebra morphism $Tridend (K)
\to \oplus_{n\geq 1} K[T_n]$ which sends the generator $x$ of $Tridend(K)$ to  $\arbreA \in T_1$. In order to prove that this
 (surjective) map is an isomorphism,
we construct explicitly its inverse in the second lemma. \hfill $\square$
\M

\N {\bf 2.7 Lemma.} {\it  The binary operations $\l , \r $ and $\cdot $ defined on $\oplus_{n\geq 1} K[T_n]$ in theorem 2.6 satisfy the
axioms of 2.1.}
\M

\N {\it Proof.} The proof is straightforward by induction on the degree. Let us show for instance that 
$$(x \l y) \l z = x \l (y * z).$$
For  $x=  x^{(0)}\vee \cdots \vee x^{(k)}$ one has
$$\eqalign{
(x \l y) \l z &=  (x^{(0)}\vee \cdots \vee (x^{(k) }* y))\l z \cr
&=  x^{(0)}\vee \cdots \vee ((x^{(k) }* y) * z). \cr
}$$
On the other hand one has
$$  x \l (y * z)= x^{(0)}\vee \cdots \vee (x^{(k) }* (y * z)).$$
Since the degree of $\x xk $ is strictly smaller than the degree of $x$, we can assume that all the relations are fulfilled for $\x xk , y $
and $z$. In particular the  associativity relation
$ (x^{(k)}* y)* z =  x^{(k)}* (y * z)$ holds.
Therefore one gets $(x \l y) \l z = x \l (y * z)$ as expected.

All the other formulas are proved similarly. \hfill $\square$
\M

\N {\bf 2.8 Lemma.} {\it Let us denote by $u$ the generator of the free dendriform trialgebra $Tridend(K)$. 
The map $\aa : \oplus_{n\geq 0} K[T_n]
\to Tridend(K)\oplus K.1$ defined inductively by
$$\aa (\vert) := 1, \qquad 
\aa (\x x0 \vee \x x1) := \aa (\x x0) \r u \l \aa(\x x1)\ ,$$
 and 
$$\aa (\x x0 \vee\cdots \vee  \x xk) := (\aa (\x x0) \r u)\cdot  \aa(\x x1 \vee \cdots \vee \x x{k-1} ) \cdot (u\l \aa (\x xk))$$
for $k\geq 2$ is a
morphism of dendriform trialgebras when restricted to $\oplus_{n\geq 1} K[T_n]$.}
\M
\N {\it Proof.}
Since it may happen  that $\x x0 =\vert $, (resp.  $\x xk =\vert $), we need to specify that $1 \r z = z = z \l 1$.
Similarly it may happen that, when $k=2$, one has $\x x1 = \vert$. So we need to specify that 
$$\aa (\x x0 \vee \vert\vee  \x x1) := (\aa (\x x0) \r u)\cdot  (u\l \aa (\x x2)).$$

For instance, one has 
$\aa (\ \arbreBC) = \aa (1 \vee 1 \vee 1) = (1\r u) \cdot (u\l 1) = u\cdot u \ $.
We want to show that 
$$\aa (x \l y) = \aa (x) \l \aa (y) , \quad \aa (x \r y) = \aa (x) \r \aa (y) ,  \quad \aa (x \cdot y) = \aa (x) \cdot \aa (y) $$
for any $x\in T_p$, $y\in T_q$, $z\in T_r$. We check the first equality, the checking of the others is similar.

Let us first check the case $ x = \x x0 \vee \x x1$. On one hand we have
$$\eqalign{
 \aa (x\l y) &= \aa (\x x0 \vee (\x x1 * y)) \cr
            &= \aa (\x x0)\r u \l \big(\aa (\x x1 * y)\big) \cr
           &= \aa (\x x0)\r u \l \big(\aa (\x x1) * \aa (y)\big), \hfill \hbox { by induction}. \cr
}$$ 
On the other hand we have
$$\eqalign{
 \aa (x) \l \aa (y) &= \aa (\x x0 \vee \x x1) \l \aa (y) \cr
                   &= \big(\aa (\x x0)\r u \l \aa ( \x x1)\big) \l \aa (y) \quad\hbox {by relation (2)},\cr
                  &= \aa (\x x0)\r u \l \big(\aa ( \x x1) *\aa (y)\big) \quad \hbox {by relation (1)}.\cr
}$$
Therefore one gets  $\aa (x \l y) = \aa (x) \l \aa (y)$ as expected.

Let us now suppose that $x = \x x0 \vee \cdots \vee \x xk$ with $k\geq 2$. One one hand we have 
$$\eqalign{
 \aa (x\l y) &= \aa ( \x x0 \vee \cdots \vee (\x xk * y)) \cr
 &= (\aa ( \x x0) \r u) \cdot \aa (\x x1 \vee \cdots \vee \x x{k-1} ) \cdot \big(u \l \aa (\x xk *y)\big) \cr
 &= (\aa ( \x x0) \r u) \cdot \aa (\x x1 \vee \cdots \vee \x x{k-1} ) \cdot \big(u \l (\aa (\x xk) *\aa (y))\big) \cr
}$$
On the other hand we have
$$\eqalign{
\aa (x)\l\aa (y) &=(\aa (\x x0 \r u)\cdot\aa (\x x1 \vee\cdots\vee\x x{k-1})\cdot (u\l\aa (\x xk))\l\aa(y)\cr
               &=(\aa(\x x0 \r u)\cdot\aa(\x x1 \vee\cdots\vee\x x{k-1})\cdot((u\l\aa (\x xk)) \l \aa(y))\cr
              &=(\aa(\x x0 \r u)\cdot\aa(\x x1 \vee\cdots\vee\x x{k-1})\cdot \big(u\l (\aa (\x xk)*\aa(y))\big)\cr
}$$
by relations (6) and (1), whence the result.

If $x= \x x0 \vee \vert \vee \x x1$, then the proof is similar and uses also the relations (6) and (1) of 2.1.
\hfill $\square$
\M

\N {\bf 2.9 Corollary.} {\it The free dendriform trialgebra  $Tridend (V)$ on the vector space $V$ is 
$$ Tridend (V) = \oplus _{n\geq 1}K[T_n]\t V^{\t n},$$ 
and the operations are induced by the operations on $ \oplus _{n\geq 1}K[T_n]$ and concatenation.} \hfill $\square$
\M
\N {\it Proof.} Follows from Theorem 2.6 by the same argument as in Corollary 1.8. \hfill $\square$
\M

\N {\bf 2.10 The family of Stasheff polytopes.} Let $\K^n$ be the Stasheff polytope (alias associahedron also  denoted $K_{n+2}$) of dimension
$n$, cf. [St]. The cells of
$\K^{n-1}$ are in one-to-one  correspondence with the planar trees with $n$ leaves. More precisely the set $T_{n,k}$
of planar trees with $n$ leaves and $n-k+1$  vertices labels the cells of dimension $k-1$ of $\K^{n-1}$. In
particular the planar binary trees are in 1-1 correspondence with the vertices. 
$$\stasheff$$
\M

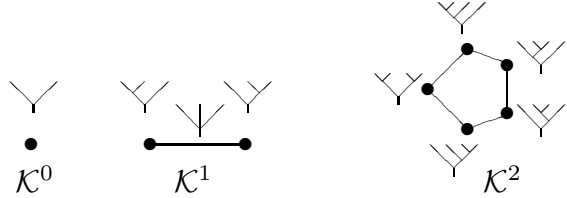
\N {\bf 2.11 The operad of dendriform trialgebras.} The operad $Tridend$ is a non-$\SS$-operad and so is completely determined by the free 
dendriform trialgebra  on one generator. The filtration on $T_n$ (cf. 2.4) is compatible with the three operations (cf. [L3, section 9.9]).
 In particular there is a functor from the category of dendriform dialgebras to the category of dendriform trialgebras 
(take $x\cdot y = 0$). The operad of dendriform dialgebras involves only the planar {\it binary} trees.

Since the operad {\it Tridend} is filtered, we can built the graded associated operad {gr\ Tridend} as follows:
$gr\ Tridend (n) := F_n Tridend / F_{n-1} Tridend$. It is clear that the 0th part of the graded operad is the operad of dendriform dialgebras. 

From the bijection between the cells of the Stasheff polytopes and the planar trees, and the description of the free dendriform trialgebra
 given in Theorem 2.6 it follows that the generating series of the family of Stasheff polytopes is equal to the generating series of the
 (filtered) operad {\it Tridend}:
$$ f^{\K }_t(x) =f^{Tridend}_t(x).$$
\M

One way of keeping track of the filtration is to introduce a type of algebra depending on a parameter $q\in K$ as follows. In the 
relations 1 and 3 of 2.1 we replace the occurences of $a \cdot b$ (where $a$ and $b$ are $x, y$ or $z$) by $q(a\cdot b)$. When  $q=1$ this is the 
dendriform trialgebra. When  $q=0$, this is (almost) the case treated by Chapoton in [Ch]. Almost because he is working with graded
 vector spaces and modify the grading for the operation $\cdot\ $, like when passing from Poisson algebras to Gerstenhaber algebras.
\M

\N {\bf 2.12 Remark.} The Stasheff polytopes form an operad (cf. [St]), which encodes the associative algebras up to
homotopy  ($A_{\infty}$-algebras). But, in this case the Stasheff polytope $\K^n$ is put in dimension $n+2$ while
for the operad $Tridend$  it is put in dimension $n+1$. In other words, in the $A_{\infty}$-algebra case a cell of the Stasheff polytope
 ${\cal K}^n$ encodes an operation on $n+2$ variables, though in the dendriform case it encodes an operation on $n+1$ variables. So they
 are completely different operads.
\BB

\N {\bf 3. Homology and Koszul duality.} In [G-K] Ginzburg and Kapranov have extended the notion of Koszul duality to binary quadratic
operads. Both operads $Trias$ and $Tridend$ are binary and quadratic, hence we can apply this theory here. In particular we can
 construct the chain complex of an associative trialgebra (resp. of a dendriform trialgebra), and also the Koszul complex of these operads. 
\M

\N {\bf 3.1 Theorem.} {\it The operad $Trias$ of associative trialgebras is dual to the operad $Tridend$ of dendriform trialgebras :
$$ Trias^! = Tridend\quad {\rm and}\quad   Tridend^! = Trias.$$}

\N {\it Proof.} Let us compute the Koszul dual of $Trias$. Since we are dealing with non-$\SS$-operads, 
that is $\P(n) = \P'(n)\t K[S_n]$ we can forget about the action of the symmetric group and work with $\P'(n)$. 
The space of generating operations is $Trias'(2) = K[P_2] = K \g \oplus K \d \oplus K \m$ . The space of operations that one 
can perform on three variables is $ K[P_2\times P_2]\oplus K[P_2\times P_2]$. This is the part of degree 3 of the free 
non-$\SS$-operad generated by $ K[P_2]$. The operad $Trias$ is completely determined by some subspace 
$R\subset K[P_2\times P_2]\oplus K[P_2\times P_2]$. Let us denote by $(\circ_1)\circ_2$ (resp. $\circ_1(\circ_2)$ ) the basis 
vectors of the first (resp. second) summand $ K[P_2\times P_2]$. Then $R$ is generated by the 11 vectors of the form 
$(\circ_1)\circ_2 - \circ_1(\circ_2)$ obtained from the 11 relations of definition of associative trialgebras (cf. 1.2). 

Let us identify the dual of $ K[P_2]$ with itself by identifying a basis vector with its dual. Then the dual operad 
$ Trias^!$ is completely determined by $R^{\perp}\subset K[P_2\times P_2]\oplus K[P_2\times P_2]$, where $R^{\perp}$ 
is the orthogonal space of $R$ under the quadratic form $\pmatrix {\Id&0\cr 0 & -\Id}$ (cf. [G-K]).

We claim that, under the identification $\g = \l$ , $\d = \r$,  $\m=\cdot $ , the space $R^{\perp}$ is the space
$R'
$  generated by the vectors obtained from the 7 relations of definition of dendriform trialgebras (cf. 2.1).
Indeed,  since $\dim K[P_2\times P_2]\oplus K[P_2\times P_2] = 18, \dim R = 11$ and $\dim R' = 7$, it is sufficient
to prove  that $\langle v, w\rangle =0$ for any basis vector $v$ of $R$ and any basis vector $w$ of $R'$. This is a
straightforward checking. We verify this equality in one case, the others are similar. Let $v=(\g)\g - \g(\g)$ which we identified
 with $(\l)\l - \l(\l)$. We get
$$\displaylines{
\langle v, (x \l y) \l z -  x \l (y * z)\rangle = 1-1 = 0  , \cr
\langle v, (x \r y) \l z - x \r (y \l z)\rangle =0  , \cr
\langle v, (x * y) \r z - x \r (y \r z)\rangle = 0 , \cr
\langle v, (x \r y)\cdot z -  x \r (y \cdot z)\rangle =0  , \cr
\langle v, (x \l y) \cdot z -  x \cdot (y \r z)\rangle =0  , \cr
\langle v, (x \cdot y) \l z -  x \cdot (y \l z)\rangle =0  , \cr 
\langle v, (x \cdot y) \cdot z - x \cdot (y \cdot z)\rangle =0  .\cr
}$$
Hence the dual of the operad $Trias$ is the operad $Tridend$.
\hfill $\square$
\M

\N {\bf 3.2 Trialgebras versus dialgebras.} There is a functor from the category of associative trialgebras to the category of
 associative dialgebras (cf. 1.4.b), that is a map from the operad  $Dias$ to the operad $Trias$. This is a map of binary quadratic
 operads. Its dual is a map from $Tridend$ to $Didend$, which gives the functor from the category of dendriform dialgebras
 to the category of dendriform trialgebras described in 2.3.b. Similarly the dual of the functor ${\bf As}\to {\bf Trias}$ is the functor 
${\bf Tridend}\to {\bf As}$.
\M

\N {\bf 3.3  Homology of associative trialgebras.} Ginzburg and Kapranov's theory of algebraic operads 
shows that there is a well-defined chain complex for any algebra $A$ over the binary quadratic operad $\P$, constructed 
out of the dual operad $\P^!$ as follows.

The chain complex of the $\P$-algebra $A$ is $C^{\P}_n(A)= \P^!(n)^*\t _{S_n} \tT An$ in dimension $n$ and the differential 
$d$ agrees, in low dimension, with the $\P$-algebra structure of $A$ $$
\cc_A(2): \P(2)\t A^{\t 2}\to A
$$
under the identification $\P^!(2)^*\cong \P(2)$.

In fact $d$ is characterized by this condition plus the fact that on the cofree $ \P^{!}$-coalgebra $C_*^{\P}(A)= \P^{!*}(A)$ it is a 
graded coderivation. 
\M

\N {\bf 3.4 Proposition.} {\it The chain complex of an associative trialgebra  $A$ is given by 
$$C^{Trias}_n(A)= K[T_n]\t A^{\t n} , \quad d = \sum_{i=1}^{i=n-1} (-1)^i d_i , $$
where $d_i(t; a_1,\cdots ,a_n) = (d_i(t); a_1,\cdots ,a_i \circ_i^t a_{i+1}, \cdots a_n),$ and $d_i(t)$ is the tree obtained from 
$t$ by deleting the $i$th leaf and where $\circ_i^t$ is given by 
$$\circ_i^t = \cases{
\g & if $ \hbox { the $i$th leaf of $t$ is left oriented} ,$\cr
 \d & if $ \hbox { the $i$th leaf of $t$ is right oriented} ,$\cr
 \m & if $ \hbox { the $i$th leaf of $t$ is a middle leaf} .$\cr
} $$}
Observe that at a given vertex of a tree there is only one left leaf, one right leaf, but there may be none or 
several middle leaves.
\M

\N {\it Proof.} First observe that this is a chain complex since the operators $d_i$ satisfy the presimplicial
relations 
$$ d_i d_j = d_{j-1} d_i \hbox { for } i<j.$$ Indeed, this relation is either immediate (when $i$ and $j$ are far apart), 
or it is a consequence of the axioms of associative trialgebras when $j=i+1$. It suffices to check the case $n=3$,
and this was done in 1.3. 

By Theorems
3.1 and 2.6 Ginzburg and Kapranov theory gives, as expected,  
$$C^{Trias}_n(A)= K[T_n]\t A^{\t n}.$$
It is clear from 1.3 that $d$ agrees with the $Trias$-algebra structure of $A$ in low dimension. Since $d$ is 
completely explicit, the coderivation property is immediate to check. \hfill $\square$ 
\M

\N {\bf 3.5 Proposition.} {\it The chain complex of a dendriform trialgebra $A$ is given by 
$$C^{Tridend}_n(A)= K[P_n]\t A^{\t n} , \quad d = \sum_{i=1}^{i=n-1} (-1)^i d_i , $$
where $d_i(X; a_1,\cdots ,a_n) = (d_i(X); a_1,\cdots ,a_i \circ_i^X a_{i+1}, \cdots a_n),$ and $d_i(X)$ is the 
image of $X$ under the map $d_i : [n] \to [n-1]$ given by $$
d_i({\u r}) = \cases{
r-1 &if $ i\leq r,$\cr
r & if $ i\geq r+1.$\cr}
$$
and where $ \circ_i^X$ is given by
$$\circ_i^X = \cases{
\cdot & if $ i-1\in X \hbox { and } i \in X,$\cr
 \succ & if $ i-1\notin X \hbox { and } i \in X,$\cr
 \prec & if $ i-1\in X \hbox { and } i \notin X,$\cr
 *	& if $ i-1\notin X \hbox { and } i \notin X.$\cr
}$$}
\M
\N {\it Proof.} Again,  observe that this is a chain complex since the operators $d_i$ satisfy the presimplicial
relations 
$$ d_i d_j = d_{j-1} d_i \hbox { for } i<j.$$ 
Indeed, this relation is either immediate (when $i$ and $j$ are far apart), 
or it is a consequence of the axioms of dendriform trialgebras when $j=i+1$. It suffices to check the case $n=3$. We actually
 do the computation in one particular case, the others are similar:
$$\displaylines{
d_1 d_2 (\{0,2\},a_1,a_2,a_3)= d_1(\{0,1\};a_1, a_2\l a_3)= (\{0\}; a_1\cdot (a_2\r a_3)\ )\ ,\cr
d_1 d_1 (\{0,2\},a_1,a_2,a_3)= d_1(\{0,1\};a_1\r a_2,a_3)= (\{0\};(a_1\l a_2)\cdot a_3\ )\ ,\cr
}$$
These two elements are equal by the fifth relation in 2.1.
 \hfill $\square$ 

\BB

\N {\bf 4. Acyclicity of the Koszul complex.} By definition the Koszul complex 
associated to the operad $Trias$ is the differential functor $Tridend \circ Trias$ from $\Vect$ to $\Vect$.
We will show that it is quasi-isomorphic to the identity functor. Equivalently we have the
\M

\N {\bf 4.1 Theorem.} {\it The homology of the free associative trialgebra  on $V$  is 
$$H^{Trias}_n(Trias(V) ) =
\cases { V & if n=1,\cr
0 & otherwise. }
$$}

\N {\bf 4.2 Corollary.} {\it The operads $Trias$ and $Tridend$ are Koszul operads.} \M
\M
\N {\bf 4.3 Corollary.} {\it The homology of the free dendriform trialgebra  on $V$  is: 
$$H^{Tridend}_n(Tridend(V) ) = \cases {V & if n=1,\cr
0 & otherwise. }
$$}

\N {\bf 4.4 Corollary.} {\it Let $f^{\K }_t(x)$ be the generating series of the Stasheff polytope (i.e. of the planar trees), 
as defined in 1.12. Then one has
$$f^{\K}_t(x) = {-(1+(2+t)x) +\sqrt { 1+2(2+t)x+t^2x^2}\over 2(1+t)x}\ .$$}
\M

\N {\it Proof of the Corollaries.} By Ginzburg and Kapranov theory [G-K] the first two Corollaries  follow from the vanishing 
of the homology of the free associative trialgebra.

The last Corollary follows from the functional equation relating the two operads and the computation of the generating series
 for the associative trialgebra operad (cf. 1.12). \hfill $\square$ 
\M

\N {\it Proof of Theorem 4.1.}
The acyclicity of the augmented complex \goodbreak $C^{Trias}_*(Trias(V) )$ is proved in several steps as follows. 
\S

1. We show that it is sufficient to treat the case $V=K$. 

2. The chain complex $C^{Trias}_*(Trias(K))$ splits into the direct sum of chain complexes $C_*(u)$,
 one for each element $u$ in $P_m$, $m\geq 1$. 

3. The chain complex $C_*(u)$ is shown to be the cell complex of a simplicial set $X(u)$. 

4. The space $X(u)$ is shown to be the join of spaces $X(v)$ for certain particular elements $v$ in $P_m$.

5. The spaces $X(v)$ are shown to be contractible by constructing a series of retractions by deformation. 

\M
\N 1. {\it First step}.
Recall from 1.8 that $Trias (V) = \bigoplus _{n\geq 1}K[P_n]\t V^{\t n}$. Therefore one has $$\eqalign{
C^{Trias}_j(Trias(V))&= K[T_j]\t \big( \bigoplus _{n\geq 1}K[P_n]\t V^{\t n} \big)^{\t j}\cr &= K[T_j]\t \bigoplus_{m\geq 1} 
(\oplus_{n_1+ \cdots +n_j =m } K[P_{n_1}\times \cdots \times P_{n_j}])\t \tT Vm .\cr} $$
Since $d$ is homogeneous in $V$, the complex $C^{Trias}_*$ splits into the direct sum of subcomplexes, 
one for each $m\geq 1$. This subcomplex is in fact of finite length and, up to tensoring by $\tT Vm$, is of the following form: 
$$\displaylines {
C_*(P_m):\qquad 0\to K[T_m\times P_1\times \cdots \times P_1] \to \cdots\hfill \cr \hfill \to \bigoplus_{n_1+ \cdots +n_j =m } 
K[T_j\times P_{n_1}\times \cdots \times P_{n_j}] \to \cdots 
\to K[T_1\times P_m].\cr
}$$
Recall that $P_1$ and $T_1$ have only one element. The case $m=1$ gives the subcomplex of length 0 reduced to $V$. This shows that 
$H_1^{Trias}(Trias(V))$ contains $V$ as expected. 

For $m\geq 2$, the differential is simply the differential of $C_*(P_m)$ tensored by the identity of $\tT Vm$, hence it is 
sufficient to prove the acyclicity of $C_*(P_m)$ to prove the theorem. 
\M
\N 2. {\it Second step}. The chain complex $C_*(P_m)$ can still be split into the direct sum of smaller complexes 
indexed by the elements $u$ of $P_m$.
Indeed, let $\alpha := (t; u_1, \cdots , u_j )\in T_j\times P_{n_1}\times \cdots \times P_{n_j}$ be a basis element. 
Under applying $j-1$ face operators successively 
to $\alpha $, we get an element $(\arbreA; u) \in T_1\times P_m$ which does not depend on the choice of the face operators because 
of the simplicial relations (cf. 3.4). Considering $t$ as an operation on $m$ variables for associative 
trialgebras, $u$ is
nothing but the result of the evaluation of $t$ on $(x,\cdots , x)$, cf. 1.3. Fixing $u$, 
let $C_*(u)$ be the chain subcomplex linearly generated by the elements $\alpha $ whose image is $(\arbreA; u) \in T_1\times P_m$. 
It is clear that $C_*(P_m)$ is the direct sum of the chain complexes $C_*(u), u\in P_m$. 

Observe that
$C_*(u)$ is of simplicial type, that is, its boundary is of the form $d = - \sum_{i=1}^{i=n-1} (-1)^i d_i$. \M

\N 3. {\it Third step}. We fix $u\in P_m$. At this point it is helpful to modify slightly 
our indexing of the faces and have them to run from 0 to $n-2$ rather than from 1 to $n-1$. We also shift the indexing of the complex 
$C_*(u)$ by 1, putting $K[T_1\times P_m]$ 
in dimension $-1$. For any generator $\alpha $ of $C_*(u)$ the faces $d_i(\alpha ), 0\leq i\leq n-2$, are still generators of $C_*(u)$. 
Hence $C_*(u)$ is the normalized augmented complex of an augmented simplicial set that we denote by $X(u)$. The nondegenerate simplices 
of $X(u)$ are the linear generators $\alpha $ of $C_*(u)$. The top dimensional ones are of the form $(t;x,\cdots, x)\in X(u)_{m-2}$ 
where $t= t_0\vee \cdots \vee t_k\in T_q$. The integer $k$ is the number of decorations (Cech signs) appearing in $u$. We denote by 
$T_{\{u\}}$ this subset of $T_m$. At the other end the augmentation set is $X(u)_{-1}= T_1\times \{u\}$ (one element). The geometric 
realization of $X(u)$ is the amalgamation of simplices $\DD^{m-2}$ 
(one for each $t\in T_{\{u\}}$) under the following rule: 

if $d_{i_k} \cdots d_{i_1}(t; x, \cdots , x)= d_{i_k} \cdots d_{i_1}(t'; x, \cdots , x)$ for some $ m-2\geq i_k\geq \cdots \geq i_1\geq 0$, then we identify the 
corresponding (oriented) faces of the simplices $t$ and $t'$. Observe that under this rule a vertex of type $i$
 is identified only with a vertex of type $i$. 
\M

\N 4. {\it Fourth step}. Let us first recall the join construction of augmented simplicial sets (cf. for instance [E-P]). 
An augmented simplicial set 
is a simplicial set $X_{.}$ together with a set $X_{-1}$ and a map $d_0: X_0 \to X_{-1}$ satisfying $d_1d_0 = d_0d_0$. The join 
of two augmented simplicial sets $X_.$ and $Y_.$ is $Z_. = X_. * Y_. $ defined by $Z_n =\bigsqcup_{p+q=n-1} X_p \times  Y_q $. 
The faces are
$$\eqalign{
d_i(x,y) &= (d_i x,y) \hbox { for } 0\leq i \leq p,\cr 
d_i(x,y) &= ( x,d_{i-p-1}y) \hbox { for } p+1\leq i \leq p+1+q,\cr 
}$$
and similarly for the degeneracies. The geometric realization of the simplicial join is the topological join 
$$X*Y = X\times I \times Y / \{(x,0,y) \sim (x',0,y), (x,1,y) \sim (x,0,y')\}.$$ In particular one has $\DD^p*\DD^q = \DD^{p+q+1}$. 

Let $u = x \cdots x \xc x \cdots x \xc x \cdots x \in P_m.$ By direct inspection we see that $X(u)$ is the simplicial 
join of the simplicial sets 
$$X(x \cdots \xc), \ X(\xc\cdots x \cdots \xc),\  \cdots ,\  X(\xc\cdots x \cdots \xc), \ X(\xc \cdots x).$$ 
The point is that
there are only one Cech signs at the extreme locations. Hence it is sufficient to show the contractibility of
$X(u)$ in  the cases $u=\xc\cdots x \cdots \xc$ and $u=\xc \cdots x$. \M

\N 5. {\it Fifth step}: the case $u=\xc\cdots x \cdots \xc$ or $u=\xc \cdots x\in P_m$. We treat in detail the case  $u=\xc \cdots x$, 
the other one is similar.

Since $u=\xc \cdots x$ the trees $t$ in $T_{\{u\}}$ are of the form

$$\figure5A$$

Hence the 0-cell $(d_0)^{m-2}(t; x, \cdots , x) = (\ \arbreC; \xc x \cdots x, x)$ is the same for all $t\in T_{\{u\}}$. We
denote this vertex by 
$P$. In other words, in the amalgamation of the $(m-2)$-simplices $(t; x, \cdots , x)$ giving $X(u)$, all the vertices of type
$m-2$ get identified to
$P$.

 We will show
that there exists a sequence of  retractions by deformation
$$X(u)= X(u)^{\langle m-2\rangle}\proj\  \cdots\  \proj X(u)^{\langle k\rangle}\ {\buildrel {\phi_k} \over \proj }\  \cdots\  \proj X(u)^{\langle
0\rangle}= P \ .$$ 
The simplicial set $X(u)^{\langle k\rangle}$ is a subsimplicial set of $X(u)$ determined by its nondegenerate $k$ simplices. It is defined
inductively as follows. We suppose that $X(u)^{\langle k\rangle}$ has been defined (the induction process begins with $k=m-2$) and we
determine $X(u)^{\langle k-1 \rangle}$. On $X(u)^{\langle k\rangle}$ we introduce the equivalence relation generated by: $\aa \sim \bb$ 
if either $d_k\aa = d_k\bb$ or 
 $d_{k-1}\aa = d_{k-1}\bb$. Then in each equivalence class we pick an element, say $\aa_0$. By definition $X(u)^{\langle k-1\rangle}$
 is made of the elements 
$d_{k-1}\aa_0 $, one for each equivalence class.

The map $\phi_k$ is defined by $\phi_k(\aa) = s_{k-1}d_{k-1}\aa_0 $.
On the geometric realization the map  $\phi_k$ consists in collapsing each $k$-simplex $\aa$ to its last face (the edge relating the vertices $k-1$ 
and $k$ collapses to a point), and
then embedding this face into $X(u)$ as $d_{k-1}\aa_0$. All the collapsing are coherent, and so assemble to give a collapsing of   $X(u)^{\langle
k\rangle}$ to  $X(u)^{\langle k-1\rangle}$, because one can verify that for each vertex of type $k-1$ in  $X(u)^{\langle k\rangle}$ there is only one edge to the
edge relating it to the vertex of type $k$, that is $P$.

Here is an illustration for  $m=4$, $u=\xc xxx$ and the planar binary trees.
$$\matrix{
&d_0 & d_1 & d_2\cr
& & & \cr
a = (\arbrequdt; x,x,x,x)& (\arbrequatre; \xc x, x, x),& (\arbrequatre; x, x\xc , x),& (\arbrequatre; x, x, x\xc ),\cr
b = (\arbreqdut; x,x,x,x)& (\arbrequatre; \xc x, x, x),& (\arbrequatre; x, \xc x, x),& (\arbrequatre; x, x, x\xc ),\cr
c = (\arbreqtud; x,x,x,x)& (\arbrequatre; \xc x, x, x),  & (\arbrecinq; x, \xc x , x),& (\arbrecinq; x, x,  \xc x),\cr
d = (\arbreqtdu; x,x,x,x)& (\arbrecinq; \xc x, x, x ),  & (\arbrecinq; x,  \xc x, x),& (\arbrecinq; x, x,  \xc x),\cr
e = (\arbrequtu; x,x,x,x)& (\arbrecinq; \xc x, x, x ),  & (\arbrecinq; x, x \xc , x),& (\arbrequatre; x, x,  \xc x).\cr
  }$$

Hence the simplices $a,b,c,d,e$ of type $\Delta^2$ are amalgamated under the following rules: 
$$ d_0(a) = d_0(b)=d_0(c),\quad d_0(d) = d_0(e), \quad d_1(c) = d_1(d), \quad d_2(a) = d_2(b).$$ 

The first two spaces of the sequence  (binary case) 
$$X(\xc xxx)= X(\xc xxx)^{\langle 2\rangle}\ \proj \ X(\xc xxx)^{\langle 1\rangle}\  \proj\ 
X(\xc xxx)^{\langle 0\rangle}= P $$ 
are shown below:
$$\courbe$$
In the planar tree case $X(u)^{\langle 2\rangle}$ is made of eleven 2-simplices,  $X(u)^{\langle 1\rangle}$ is made of seven 1-simplices and 
 $X(u)^{\langle 0\rangle}$ is made of one 0-simplex (namely $P$).
 
Since each map $\phi_k$ is a retraction by deformation, the space $X(u)= X(u)^{\langle m-2\rangle}$ has the same homotopy type as $
X(u)^{\langle 0\rangle}= P$ hence it is contractible. \hfill $\square$
\BB

\N {\bf 5. Cubical trialgebras and hypercubes.} 
\M

One can also associate a type of trialgebras to the family of hypercubes. 
Once the correct relations are found the proof follows the same pattern as in the previous sections. It turns out that the 
associated operad is self-dual, so the generating series, which is $f^{I}_t(x) = {-x \over 1+(2+t)x}\ ,$ is its own
inverse, a fact which is immediate to check:
$ f^{I}_t(f^{I}_t(x) ) = x\ .$
\M

\N {\bf 5.1 Definition.} A {\it cubical trialgebra} is a vector space $A$ equipped with 3 binary operations : 
$\g$ called {\it left}, $\d$ called {\it right} and $\m$ called {\it middle}, satisfying the following 9 relations :
$$(x\circ_1 y)\circ_2 z = x\circ_1 (y\circ_2 z)$$ where $\circ_1$ and $\circ_2$ are either $\g$ or $\d$ or $\m$. \M

We obtain the definition of a {\it cubical dialgebra} by restricting ourself to the first two operations 
(this structure has been considered earlier by B. Richter [Ri]). We denote by {\bf Tricub} and {\bf Dicub} the associated 
categories of algebras.
There is an obvious functor
$${\bf As} \to {\bf Tricub}$$
consisting in putting $x\g y = x\d y =x\m y = x y $.
\M
Let $Q_n$ be the set of cells of the hypercube $I^n$, where $I$ is the interval $[-1,1]$. Alternatively $Q_n$ can be described 
as $\{-1, 0, +1\}^n$ or $\{\g, \m, \d\}^n$. Obviously $Q_n$ is graded by the 
dimension of the cells (resp. the numbers of 0's or $\m$ signs). 
\M

\N {\bf 5.2 Proposition. } {\it The free cubical trialgebra on one generator, 
$Tricub(K)= \bigoplus_{n\geq 1}Tricub(n)$ is such that $Tricub(n)= K[Q_{n-1}]$ with operations: 
$$\eqalign{
a\g b &= (a, -1, b)\in Q_{p+q-1}, \cr
a\m b &= (a, 0, b)\in Q_{p+q-1}, \cr
a\d b &= (a, +1, b)\in Q_{p+q-1}, \cr
}$$
for $a\in Q_{p-1}$ and $b\in Q_{q-1}.$}
 \M
\N {\bf 5.3 Theorem.} {\it The operad $Tricub$ is self-dual.}
 \M
(Observe that $18-9 = 9$.)\hfill $\square$ 
\M

\N {\bf 5.4 Cubical trialgebras and associative algebras.} By Koszul duality the functor ${\bf As} \to {\bf
Tricub}$ gives a functor ${\bf Tricub} \to {\bf As}$ since both operads are self-dual. It is immediately seen that
it is given by putting
$ x*y := x\g y + x\d y + x\m y$.  So a cubical trialgebra is an associative algebra for which the associative
operation is the sum of  three operations and the associative relation splits into 9 relations.
\M

\N {\bf 5.5 Proposition. } {\it The homology of a cubical trialgebra $A$ is given by the following chain complex
$C^{Tricub}_n(A)$: 
$$C^{Tricub}_n(A)= K[Q_n]\t A^{\t n} , \quad d = - \sum_{i=1}^{i=n-1} (-1)^i d_i , $$
where $d_i(X; a_1,\cdots ,a_n) = (d_i(X); a_1,\cdots ,a_i \circ_i^X a_{i+1}, \cdots a_n),$ and the element $d_i(X)$ is obtained from 
$X$ by deleting the $i$th coordinate $X_i$, and the operation $ \circ_i^X$ is given by $$
\circ_i^X = \cases{
\g & if $ X_i = -1,$\cr
\m & if $ X_i = 0,$\cr
\d & if $ X_i = +1.$\cr}
$$}
\N {\bf 5.6 Theorem.} {\it Let $Tricub(V)$ be the free cubical trialgebra on $V$. Its homology is 
$$H^{Tricub}_n(Trias(V) ) = \cases {
V & if n=1,\cr
0 & otherwise. }
$$}

\N {\bf 5.7 Corollary.} {\it The operad $Tricub$ is a Koszul operad.} 
\M
\N {\it Proof.} The same arguments as in the proof of theorem 4.1 lead to the chain complex 
$$\displaylines {
0\to K[Q_m\times Q_1\times \cdots \times Q_1] \to \cdots\hfill \cr
 \hfill \to \bigoplus_{n_1+ \cdots +n_j =m } K[Q_j\times Q_{n_1}\times \cdots \times Q_{n_j}] \to \cdots 
\to K[Q_1\times Q_m].\cr
}$$
This complex is the direct sum of complexes $C_*(u)$, one for each generator $u$ of $Q_m$. By direct inspection we see that 
$C_*(u)$ is nothing but the normalized augmented chain complex of the standard simplex $\Delta ^{m-1}$, hence it is acyclic. 

\hfill $\square$ 

\BB

\centerline {\bf References}
\B\N [A] Aguiar, M. {\it Pre-Poisson algebras}. Lett. Math. Phys. 54 (2000), no. 4, 263--277.
\S

\N [Ca] Cartier, P. {\it Fonctions polylogarithmes, nombres polyzetas et groupes pro-unipotents.} S\'eminaire Bourbaki, mars 2001, 
expos\'e 885.
\S
\N [Ch] Chapoton, F., {\it Construction de certaines op\'erades et big\`ebres associ\'ees aux polytopes de Stasheff et hypercubes},
    Trans. A.M.S 354  ({\oldstyle 2002}),  63--74.
\S
\N [E-P] Ehlers, P. J.; Porter, T. {\it Joins for (augmented) simplicial sets}. J. Pure Appl. Algebra 145 (2000), no. 1, 37--44. 
\S
\N [G-K] Ginzburg, V. and  Kapranov, M. {\it Koszul duality for operads}. Duke Math. J. 76 (1994), no. 1,
203--272.
\S
\N [H] Hoffman, M. E. {\it
Quasi-shuffle products.}
J. Algebraic Combin. 11 (2000), no. 1, 49--68. 
\S
\N [L1] Loday, J.-L. {\it Alg\`ebres ayant deux op\'erations associatives (dig\`ebres)}. C. R. Acad. Sci.
Paris S\'er. I Math. 321 (1995), no. 2, 141--146. 
\S
\N [L2] Loday, J.-L. {\it Dialgebras}, in ``Dialgebras and related operads'', Springer Lecture Notes in Math., vol
1763 (2001), 7--66.
\S
\N [L3] Loday, J.-L. {\it Arithmetree}, J. Algebra (to appear).
\S
\N [LR1] Loday, J.-L.; Ronco, M. O. {\it Hopf algebra of the planar binary trees}. Adv. Math. 139 (1998),
no. 2, 293--309.
\S
\N [LR2] Loday, J.-L.; Ronco, M. O. {\it Une dualit\'e entre simplexes standards et polytopes de Stasheff},
 C.R.Acad.Sci. Paris t. 333, S\'er. I (2001), 81 -86.
\S
\N [O] Operads: Proceedings of
Renaissance Conferences (Hartford, CT/Luminy, 1995), 37--52, Contemp. Math., 202, Amer. Math. Soc., Providence, RI, 1997. 
\S
\N [R] Ronco, M. {\it A Milnor-Moore theorem for dendriform Hopf algebras}. C. R. Acad. Sci. Paris S\'er. I Math. 332 (2000), 109-114. 
\S
\N [R] Ronco, M. {\it Eulerian idempotents and Milnor-Moore theorem for certain noncocommutative Hopf algebras}. J. Algebra (to appear).
\S
\N [Ri] Richter, B. {\it Dialgebren, Doppelagebren und ihre Homologie}. Diplomarbeit, Bonn Universit\" at, unpublished. 
\S
\N [St] Stasheff, James Dillon. {\it Homotopy associativity of $H$-spaces}. I, II. Trans. Amer. Math. Soc. 108 (1963), 275-292;
 ibid. 293--312. 
\B
\N JLL : Institut de Recherche Math\'ematique Avanc\'ee, 

CNRS et Universit\'e Louis Pasteur, 7 rue R. Descartes, 

67084 Strasbourg Cedex, France

E-mail : loday@math.u-strasbg.fr
\M

\N MOR : Departamento de Matem\'atica

Ciclo B\'asico Com\'un, Universidad de Buenos Aires 

Pab. 3 Ciudad Universitaria Nu\~nez

(1428) Buenos-Aires, Argentina

E-mail : mronco@mate.dm.uba.ar
\end